\documentclass[10pt,leqno,twoside]{article} 
\usepackage{amsmath,fancyhdr,amssymb,bm} %
\def\txtwidth#1{\setlength\textwidth{#1truein}}   
\def\ooffset#1{\setlength\oddsidemargin{#1truein}}   
\def\eoffset#1{\setlength\evensidemargin{#1truein}}   
\def\tpmargin#1{\setlength\topmargin{#1truein}}   
\def\tpskip#1{\setlength\topskip{#1truein}}   
\def\txtheight#1{\setlength\textheight{#1truein}}   
\def\ftskip#1{\setlength\footskip{#1truein}}   
\def\ha{{\textstyle\frac{1}{2}}}
\def\R{\mathbb{R}}  
\def\C{\mathbb{C}}   \def\ha{\sfr 1/2}
\def\chir{\raise2pt\hbox{$\chi$}}

\DeclareMathSymbol{\pri}{\mathcal}{symbols}{"30}
\renewcommand{\prime}{\hskip0.8pt\pri\hskip-0.4pt}

\renewcommand{\thefootnote}{\fnsymbol{footnote}}
\def\res{\hbox{ {\vrule height .22cm}{\leaders\hrule\hskip.2cm}}\hskip5.0\mu}

\def\op#1{\operatorname{\text{\rm #1}}}

\newcommand{\prelemskip}{\vskip 5pt plus 1.5pt minus 1.5pt}
\newcommand{\postlemskip}{\vskip 4.5pt plus 1pt minus 0.5pt} 
\newenvironment{state}[1]{\par\prelemskip\noindent {\bf#1} % 
\textit\bgroup\abovedisplayskip=5pt plus 3pt minus 1pt \belowdisplayskip=5pt plus 3pt minus 1pt}%  
{\egroup \par\ifdim\lastskip<\medskipamount \removelastskip\penalty55\postlemskip\fi}

\def\op#1{\operatorname{\text{\rm #1}}}

\def\tint{\text{\small$\int\hskip-2pt$}}  
\newcommand{\sing}{\op{sing}}

\newcommand{\graph}{\op{graph}}
\newcommand{\dist}{\op{dist}}

\newcommand{\dvg}{\op{div}}
\newcommand{\support}{\op{support}}

\tpmargin{-0.35}\tpskip{0.45}\txtheight{8.7}\ftskip{0.2}   %vertical dim (ordered in truein) 
\txtwidth{5.2}\ooffset{0.65}\eoffset{0.65}                            %horizontal dimensions 
\def\R{\mathbb{R}} \def\C{\mathbb{C}}   
\def\ha{{\textstyle\frac{1}{2}}}

\newcommand\wtilde[1]{\widetilde #1}

\title{\vskip-.2in A general asymptotic decay lemma \\ for elliptic problems} \date{} \author{{\scshape Leon
Simon}\thanks{This paper is dedicated to S.-T.~Yau on the occasion of his 60'th birthday, in recognition of his
many contributions to the development of the field of geometric analysis, and in sincere appreciation of 
a  friendship which has spanned several decades.
The present work is a revision of an earlier (unpublished) preprint and has been supported by NSF grants
DMS--0406209 \& DMS--0104049 at Stanford University}}%

\begin{document}

\maketitle 

\vspace{-.4in}

\section*{Introduction} 

% draft of intro at end.

Asymptotic decay and growth theorems are fundamental in the study of geometric variational problems.  For
example in the study of minimal surfaces the pioneering work of De~Giorgi, Reifenberg, Almgren and Allard
depended on proving appropriate asymptotic decay lemmas near ``regular points.''  In later work, asymptotic
behavior near singularities has proved to be a key ingredient in attempts to understand the nature of the singular
set.

While much has been achieved, nevertheless many basic questions concerning asymptotics near singularities
remain open.  For example, perhaps the most famous and basic of all open questions concerning asymptotics,
there is the question of existence of a unique tangent cone for a minimal surface at each of its singular
points---that is, the question of whether a singular minimal surface (or more generally a stationary integral
varifold) is asymptotic to a cone on approach to each of its singular points.

We make no attempt here to give a systematic survey of the various works which address such questions, some
references for which would include for example \cite{Rei60}, \cite{Deg61}, \cite{Alm66}, \cite{Alm00},
\cite{BomDG69}, \cite{All72}, \cite{BomG72}, \cite{All75}, \cite{SchSY75}, \cite{Tay76}, \cite{Mir77},
\cite{HarS79}, \cite{SchS81}, \cite{Giu83}, \cite{Whi83}, \cite{Sim83a}, \cite{Sim87a}, \cite{AdaS88},
\cite{ChaS88}, \cite{Sim89}, \cite{Whi92}, \cite{Sim93}, \cite{Sim95a}, \cite{Sim95b}.  Rather here we will
discuss one general, but technically reasonably straightforward, asymptotically decay lemma, in the hope that it
will provide part (albeit a small part) of an effective introduction to the more technical works mentioned above.

The general asymptotic decay/growth theorem discussed here is applicable to various geometric variational
problems, and gives general criteria for establishing growth and decay properties in the presence of singularities. 
The main results (Theorems~1,2 in~\S1) can be applied to positive supersolutions $u$ of equations of the form
$\Delta_{\! M}u+b\cdot \nabla u+(q+a)u=0$ with $q\ge 0$ and $a,b$ ``small'' perturbation terms, provided that
the submanifold $M$ is part of a suitable ``regular multiplicity~1 class'' of submanifolds and, in the case of
Theorem~2, provided that the pair $M,q$ is ``asymptotically conic'' in the appropriate sense.  The terminology
is made precise in~\S1 below.

One of the principal technical ingredients is the partial Harnack theory developed in~\S5, which is key to
ensuring that ``concentration of $L^{p}$-norm'' does not occur. The main theorem (Theorem~1) is stated in \S1
and proved in~\S6.

The applicability of the main theorem to interesting geometric problems is illustrated in \S7, where we describe
how the general theorem applies to give growth estimates for entire and exterior solutions of the minimal
surface equation---i.e.\ for solutions of the minimal surface equation which are either $C^{2}$ on all of $\R^{n}$
or else $C^{2}$ on $\R^{n}\setminus\Omega$, where $\Omega$ is a bounded open subset of $\R^{n}$, in case
the gradient is unbounded.  (If the gradient of an exterior solution is bounded then it has a limit at infinity by a
result of Bers~\cite{Ber51} for $n=2$ and by~\cite{Sim87a} for $n\ge 3$; of course as pointed out
in~\cite{BomDM69}, entire solutions of bounded gradient are actually linear by the $C^{1,\alpha}$ estimate for
solutions of quasilinear elliptic equations (\cite[Th.~13.1]{GilT83}), which gives
$R^{\alpha}[Du]_{\alpha,B_{R}}\le C$, $C=C(n)$, whence by letting $R\to \infty$ we obtain
$[Du]_{\alpha,\R^{n}}=0$, i.e.\ $Du$ is constant on $\R^{n}$.)  The result obtained in~\S7 is summarized in the
following theorem:

\begin{state}{\bf{}Theorem.} Suppose $\Omega$ is a bounded open subset of\/ $\R^{n}$ and $u$ is a $ C^{2}$
solution of the minimal surface equation on $\R^{n}\setminus\Omega$ such that $|Du|$ is not bounded.  Then
for each $\gamma<\gamma_{0}\equiv{n-3\over{}2}-\sqrt{({n-3\over{}2})^{2}-(n-2)}$ there are constants
$C,R_{0}>0$ (depending on $u$) such that
\begin{align*}%
  R^{-n}\int_{S_{R}}\nu_{n+1}\, d{\cal{}H}^{n} &\le C R^{-\gamma}, \quad \forall\, R\ge R_{0}\\ %
  R^{-n}\int_{S_{R}}|Du|\, d{\cal{}H}^{n}         &\ge{}C R^{\gamma},\quad\forall\, R\ge R_{0}.%
\end{align*}%
Here $S_{R}=\{(x,u(x))\in(\R^{n}\backslash\Omega)\times\R:|x|^{2}+u^{2}(x)<R^{2}\} \equiv{}G\cap
\{(x,y)\in\R^{n+1}:|(x,y)|<R\}$, with $G=\{(x,u(x)):x\in\R^{n}\setminus\Omega\}=\graph{}u$,
$\nu_{n+1}=(1+|Du|^{2})^{-1/2}$ is the $(n+1)$'st component of the upward pointing unit normal
$\nu=(1+|Du|^{2})^{-1/2}(-Du,1)$ of the graph $G$ (viewed as the restriction to $G$ of a function of
$(x,y)\in\Omega\times\R$ which is independent of the variable $y$), and ${\cal{}H}^{n}$ is
$n$-dimensional Hausdorff measure on $G$.
\end{state}
 
We prove the first inequality above in \S7 as a consequence of the main decay estimate in Theorem~2 of \S1
below. The second inequality is a consequence of the first by virtue of the Cauchy-Schwarz inequality and the
fact that there is a fixed constant $C=C(n)$ such that $C^{-1}R^{n}\le |S_{R}|\le CR^{n}$ for all
$R>2\,\text{diam}\,\Omega$.

The above theorem extends work of Caffarelli, Nirenberg, and Spruck \cite{CafNS90}, Ecker \& Huisken
\cite{Eck90}, and Nitsche \cite{Nit89} with respect to the relevant growth exponent and also with respect to the
information it gives with regard to the generality of the set of points $x$ where an inequality like $|Du(x)|\ge C
R^{\gamma}$ must hold. In particular the exponent $\gamma_{0}$ in the above theorem is best possible in
general because the original examples of non-linear entire solutions of the MSE constructed by~\cite{BomDG69}
(see also the discussion of~\cite{Sim89}) have exactly this growth. For further discussion we refer to
\cite{Sim08a}.

Other applications of the main theorem here will be described elsewhere---see in particular \cite{Sim08b}.

\section[Scale invariant compact classes of submanifolds] {Main Results}

Let ${\cal{}P}$ be a collection of properly embedded $C^{1}$ submanifolds $P$ in $\R^{N}$ and corresponding to
each $P\in {\cal{}P}$ we assume there is given an open subset $U_{P}$ of $\R^{N}$ which contains $P$.  The
collection ${\cal{}P}$ will be called a \textit{regular multiplicity~1 class} if the following conditions are
satisfied, in which we use the notation that $\breve B_{\rho}(y)$ is the open ball in $\R^{N}$ (the notation
$B_{\rho}(y)$ being reserved for the closed ball):

\medskip

{\bf 1.1 (Reducibility of ${\cal{}P}$):} $P\in {\cal{}P}$ and $\breve B_{\rho}(y)\subset U_{P}$ with
$P\cap\breve{B}_{\rho}(y)\neq \emptyset\Longrightarrow$ each connected component of $P\cap \breve
B_{\rho}(y)$ is also in ${\cal{}P}$ with $U_{P\cap \breve B_{\rho}(y)}= \breve B_{\rho}(y)$.

{\bf 1.2 (Scale invariance of ${\cal{}P}$):} $P\in{\cal{}P}\Longrightarrow\eta_{y,\rho}(P)\in{\cal{}P}$ for each
$y\in \R^{N}$ and $\rho>0$, and $U_{\eta_{y,\rho}P}=\eta_{y,\rho}(U_{P})$; here $\eta_{y,\lambda} :
\R^{N}\to \R^{N}$ is defined by $\eta_{y,\lambda}(x)= \lambda^{-1}(x-y)$.

{\bf 1.3 (Regularity property of ${\cal{}P}$):} ${\cal{}H}^{n-2}(\sing P\cap K)<\infty\,\forall \mbox{ compact
}K\subset U_{P} \mbox{ and } P\in {\cal{}P}$, where $\sing P$ is the singular set of $P$ defined by
$\sing P= U_{P}\cap \overline P\setminus P$, where $\overline P$ is the closure of $P$ as a subset of $\R^{N}$.
% We also assume that there are no removable singularities, so that XXX don't need this %

{\bf 1.4 (Compactness of ${\cal{}P}$):} $\sup_{\{P\in {\cal{}P}:\, U_{P}\supset \breve B_{1}\}}{\cal{}H}^{n}(P\cap
B_{\theta})<\infty\,\forall\,\theta\in (0,1)$, and for every sequence $\{P_{k}\}\subset {\cal{}P}$ with
$U_{P_{k}}\supset\breve{B}_{1}$ for each $k$, there is a subsequence of $\{P_{k}\}$ converging locally in
$\breve B_{1}$ in the Hausdorff distance sense to either the empty set or to some $P\in {\cal{}P}$ with
$U_{P}\supset\breve B_{1}$, and in the latter case we also require that locally, in a neighborhood of each
compact subset $K\subset P\cap \breve B_{1}$, the convergence holds in the $C^{1}$-sense that there is a fixed
open set $U$ in $\R^{N}$ with $K\subset U$ and a sequence $\Psi_{k}$ of $C^{1}$ diffeomorphisms $U\to U$
with $\Psi_{k}$ converging to the identity map on $U$ in the $C^{1}$ norm and with $\Psi_{k}(P\cap
U)=P_{k}\cap U$ for each sufficiently large $k$.

\medskip

{\bf{}Remark:} Notice that the above enables us to make good sense of statements like $f_{k}\to f$ locally in
$L^{p}$ or locally in $C^{0}$ on $P$, even if the $f_{k}$ are actually defined on $P_{k}$ (with $P_{k}\to P$ as
in~1.4) rather than on the fixed domain $P$. For example $f_{k}\to f$ locally in $C^{0}$ means that
for each compact $K\subset P$ we have $f_{k}\circ \Psi_{k}|K$ converges uniformly to $f|K$,
where $\Psi_{k}$ are as in~1.4.

\bigskip

${\cal{}C}$ will denote the set of \emph{cones}\/ $\C$ in ${\cal{}P}$, so that $\C\in{\cal{}C}$ means
$U_{\C}=\R^{N}\setminus \{0\}$ and $ \eta_{0,\rho}\C=\C\,\,\forall\, \rho>0$, where, here and subsequently (as
in~1.2), $\eta_{y,\rho}$ denotes the translation and scaling given by
$$%
\eta_{y,\rho}(x)=\rho^{-1}(x-y).   
$$%
We also let ${\cal{}E}$ be the corresponding class of
$(n-1)$-dimensional submanifolds of $S^{N-1}$:  
 $$% 
{\cal{}E}=\{\Sigma=\C\cap S^{N-1}:\C\in {\cal{}C}\},
 $$%
 equipped with the Hausdorff distance metric $d$. Evidently, in view of~1.4, ${\cal{}E}$ is then a compact metric
space.  Subsequently we let
$$%
{\cal{}C}_{0} \text{ and } {\cal{}E}_{0} \text{ denote compact subsets of } {\cal{}C},\,{\cal{}E} \text{ respectively}
\leqno{1.5} 
$$%
and correspondingly a collection 
$$%
\left\{\begin{aligned}%
&{\cal{}Q}_{0}=\{q_{\Sigma}\}_{\Sigma\in {\cal{}E}_{0}} \text{ with } %
 q_{\Sigma}\ge 0,\,\,\, q_{\Sigma}\text{ locally bounded, %
                                                              measurable on }\Sigma,   \\ %
&\Sigma_{k},\Sigma\in{\cal{}E}_{0} \text{ with } \Sigma_{k}\to \Sigma   %
                                            \text{ with respect to the Hausdorff distance } \\ %
\noalign{\vskip-2.5pt}
&\hskip0.8in \text{metric } d \Rightarrow  q_{\Sigma_{k}}\to q_{\Sigma}  %
                                                \text{ uniformly on compact subsets of $\Sigma$.} %
\end{aligned}\right.%
\leqno{1.6}
$$%
(i.e., the $q_{\Sigma}$ depend locally uniformly on $\Sigma$ with respect to the Hausdorff distance metric on
${\cal{}E}_{0}$). For $\Sigma\in{\cal{}E}_{0}$ and $q_{\Sigma}\in{\cal{}Q}_{0}$ as above, for each connected
component $\Sigma_{\ast}$ of $\Sigma$ we let $\lambda_{1}(\Sigma_{\ast})$ be the ``minimum eigenvalue'' of
the operator $-(\Delta_{\Sigma}+ q_{\Sigma})$ on the component $\Sigma_{\ast}$:
 $$%
 \lambda_{1}(\Sigma_{\ast})=  %
   \inf_{\zeta\in C^{\infty}_{c}(\Sigma_{\ast}),\,\|\zeta\|_{L^{2}(\Sigma_{\ast})}=1} %
                                   \int_{\Sigma_{\ast}}\Bigl(|\nabla\zeta|^{2}-q_{\Sigma}\zeta^{2}\Bigr). %
\leqno{1.7}
 $$%
 The reader should note that perhaps the word ``eigenvalue'' is misleading here since although the real number
$\lambda_{1}(\Sigma_{\ast})$ exists, there may be no $\varphi\in W^{1,2}(\Sigma_{\ast})$ with
$-(\Delta_{\Sigma_{\ast}}\varphi+ q_{\Sigma_{\ast}})\varphi=\lambda_{1}(\Sigma_{\ast})\varphi$, even weakly, on
$\Sigma_{\ast}$.  Of course if $\Sigma_{\ast}$ is a compact smooth manifold (i.e.\ $\sing
\Sigma_{\ast}=\emptyset$) then the usual Hilbert space applied in the space $W^{1,2}(\Sigma_{\ast})$ guarantees
such a function $\varphi$ does indeed exist, and in this case by elliptic regularity theory
(\cite[\S8.8--\S8.10]{GilT83}) it will be continuous and everywhere non-zero on $\Sigma_{\ast}$.  In general,
when $\sing \Sigma\neq \emptyset$, the De~Giorgi Nash Moser theory does guarantee the existence of a positive
$\varphi\in W^{1,2}_{\text{loc}}(\Sigma_{\ast})\cap L^{p}(\Sigma_{\ast})$ solution of the equation for
$p<{n\over{}n-2}$, as we discuss below.

With $\lambda_{1}(\Sigma_{\ast})$ as in~1.7, we define
 $$%
   \lambda_{1}(\Sigma)= \max\{\lambda_{1}(\Sigma_{\ast}):\Sigma_{\ast} %
   \text{ is a connected component of } \Sigma\}  %
 $$%
(notice that this makes sense, because, as we show in~2.4 of the next section, there are only finitely many
connected components $\Sigma_{\ast}$ of $\Sigma$), and we let
$$%
\lambda_{1}({\cal{}E}_{0}) = \sup_{\Sigma\in{\cal{}E}_{0}}\lambda_{1}(\Sigma). 
 \leqno{1.8}
$$%
The main theorems below relate to asymptotics for positive supersolutions $u$ of suitable elliptic equations on
various subdomains of $M\in{\cal{}P}$. Specifically, we assume $\tau\in (0,{1\over{}4}]$ (to be specified in the
main theorem) and $U_{M}\supset B_{3/2}\setminus \breve B_{\tau}$, and the main theorem (Theorem~1)
below assumes $u$ is given on $M\cap \breve B_{3/2}\setminus B_{\tau}$ with
 $$%
 u\in{}W^{1,2}_{\text{loc}}(M\cap \breve B_{3/2}\setminus B_{\tau})\setminus\{0\},\,\,u\ge{}0 \text{ a.e.},\,%
                                                                    \Delta_{\!  M}u+b\cdot\nabla_{\! M}u +(q +a )u\le{}0%
 \leqno{1.9}
 $$%    
 on $M\cap \breve B_{3/2}\setminus B_{\tau}$, where $q\in L^{\infty}_{\text{loc}}(M\cap \breve
B_{3/2}\setminus B_{\tau})$ with $q\ge 0$ and where $a:M\cap \breve B_{3/2}\setminus B_{\tau}\to \R$ and
$b:M\cap \breve B_{3/2}\setminus B_{\tau}\to \R^{n}$ are given locally bounded measurable functions. Notice
that since $\sing M\equiv \overline M\setminus M$ is in general non-empty, the fact that $q\in
L^{\infty}_{\text{loc}}$ of course leaves open the possibility that $q$ can be unbounded in the neighborhoods
$\sing M$.  Of course the inequality in~1.9 is to be interpreted in the weak sense that
 $$%
 \int_{\! M}\Bigl(-\nabla_{\! M} u\cdot \nabla_{\! M} \zeta+ %
 b\cdot \nabla_{\! M}u \zeta + (q+a)\,\zeta u\Bigr) \le 0, %
 \quad \zeta\in C^{1}_{c}(M\cap \breve B_{3/2}\setminus B_{\tau}),\,\,\zeta\ge 0.  %
 $$%
 Here $\zeta\in C^{1}_{c}(M\cap \breve B_{3/2}\setminus B_{\tau})$ while $\sing M=\overline M\setminus M$,
so $\support\zeta\cap\sing{}M=\emptyset$, hence, in the above, and subsequently, there is no a-priori assumption
on how $u$ behaves on approach to $\sing M$.

The functions $a,b$ should here be though of as ``perturbation terms'' and are included for reasons of
generality. Such terms are not needed (and can be taken to be identically zero) in the  application to
solutions of the minimal surface equation discussed in~\S7.  We shall in any case for the main theorems
(Theorems 1--3) need to assume $a,b$ small; we quantify this below.

The main growth theorem below considers the case when $M$ is close to a cone $\C\in {\cal{}C}_{0}$ in an
annular region $\breve B_{3/2}\setminus B_{\tau}$ in the sense
$$%
d(M\cap \breve B_{3/2}\setminus B_{\tau},\C\cap \breve B_{3/2}\setminus B_{\tau})< \tau %
 \leqno{1.10}
 $$%
 where $d$ is the Hausdorff distance metric for subsets of $\R^{N}$.  With $\Sigma=\C\cap S^{N-1}$ the
corresponding submanifold in ${\cal{}E}_{0}$, we also need to assume that the perturbation terms $a,b$ are
suitably small and that the function $r^{2}q$ is close, on $M\cap \breve B_{3/2}\setminus B_{\tau}$, to the
corresponding $q_{\Sigma}\in{\cal{}Q}_{0}$ of~1.6 in the sense that
$$%
\left\{\begin{aligned}%
&\sup_{\! \{x\in\C :\dist(x,\sing\C)\ge\tau\}\cap  B_{3/2}\setminus \breve B_{\tau}} \hskip-.2in  %
  \bigl(r^{2}(|a|^{1/2}+|b|)\circ \Psi+|r^{2}q\circ \Psi-q_{\Sigma}|\bigr) \le \tau \\ % 
& \hskip.75in \| |a|^{1/2}+|b|\|_{L^{n+\alpha}(M\cap \breve B_{3/2}\setminus B_{\tau})} \le \beta, %
\end{aligned}\right.%
\leqno{1.11} %
$$%
where $\Psi:U\to U$ is a $C^{1}$ diffeomorphism of some open $U$ containing the compact set $K=\{x\in\C
:\dist(x,\sing\C)\ge \tau\}\cap B_{3/2}\setminus \breve B_{\tau}$ with
$$%
\Psi(K)=\{x\in M :\dist(x,\sing\C)\ge \tau\}\cap B_{3/2}\setminus \breve B_{\tau},\text{ and }
\|\Psi-I\|_{C^{1}(U)}<\tau.
$$%
In all that follows, $L^{p}(\Omega)$ norms (with $\Omega\subset M$) always denote the normalized
$L^{p}$-norm, with normalizing factor chosen so that the indicator function of $\Omega$ has norm~$1$; thus
 $$%
 \|f\|_{L^{p}(\Omega)} =\Bigl(|\Omega|^{-1}\int_{\Omega}|f\,|^{p}\Bigr)^{1/p},\quad
|\Omega|={\cal{}H}^{n}(\Omega).
 $$%
We are now ready to state the main growth theorem. In the statement,
$$%
\gamma_{0}=\begin{cases}%
  {{n-2}\over{}2}-\sqrt{({n-2\over{}2})^{2}+\lambda_{1}({\cal{}E}_{0})} & \text{ if }  %
                                                 \lambda_{1}({\cal{}E}_{0})   \ge -(\tfrac{n-2}{2})^{2}\\    %
\tfrac{n-2}{2} & \text{ otherwise,} %
 \end{cases} %
$$%
where $\lambda_{1}({\cal{}E}_{0})$ is as in~1.8. 

 \begin{state}{\bf Theorem 1 (Main Growth Theorem.)} %
   For each $\alpha,\beta>0$, $p\in [1,{n\over{}n-2})$, and $\gamma<\gamma_{0}$ ($\gamma_{0}$ as above),
there is $\tau=\tau(\gamma,p,{\cal{}C}_{0},\alpha,\beta)\in (0,{1\over{}8})$ and $\rho=\rho(\gamma,p,
{\cal{}C}_{0},\alpha,\beta)\in (2\tau,{1\over{}4})$ such that if $M\in {\cal{}P}$, $\C\in{\cal{}C}_{0}$, and if
1.1--1.11 all hold, then
 $$%
 \|u\|_{L^{p}(M\cap \breve B_{\rho}\setminus B_{\rho/2})} \ge%
                \rho^{-\gamma} \|u\|_{L^{p}(M\cap \breve B_{1}\setminus B_{1/2})}.%
 $$%
\end{state}

{\bf Remarks: (1)} A key point here is that the constants $\tau,\rho$ do not depend on the particular $M,\C$
under consideration, so Theorem~1 can be applied uniformly across a large family of different $M$ and $\C$;
this will be used in the proof of the corollaries below. Of course the theorem still has content in the special case
when ${\cal{}E}_{0}$ consists of just one element $\Sigma\in {\cal{}E}$, and in this case we have
$\lambda_{1}({\cal{}E}_{0})=\lambda_{1}(\Sigma)$.

\smallskip

{\bf{}(2)} We shall show in~4.3 below that in fact under the hypotheses 1.9, 1.10, 1.11 we always automatically
have a lower bound $\lambda_{1}({\cal{}E}_{0})\ge -\bigl({n-2\over{}n}\bigr)^{2}-\epsilon(\tau)$, with
$\epsilon(\tau)\downarrow 0$ as $\tau\downarrow 0$, and of course trivially $\lambda_{1}(\Sigma)\le 0$
because $q_{\Sigma}\ge 0$, so the constant $\gamma_{0}$ in the above theorem is a well-defined real number in
the interval $[0,{n-2\over{}2}]$ and in fact for $n\ge 3$ $\gamma_{0}>0$ unless $\lambda_{1}(\Sigma)=0$, which
evidently occurs only when $q_{\Sigma}=0$ a.e.\ on $\Sigma$.

\bigskip

We shall give the proof of Theorem~1 in~\S6, after the necessary preliminaries are established. For the moment
we establish a corollary of Theorem~1.  This corollary applies to $(M,a,b,q)$ which are ``asymptotically conic''
either at 0 or $\infty$ in the following sense:

\medskip

{\bf 1.12 Definition:} For $M\in {\cal{}P}$, $q:M\to [0,\infty)$, $a:M\to \R$, $b:M\to \R^{n}$:

(a) $(M, q,a,b)$ is asymptotically conic at $x_{0}\in \overline M\cap U_{\!  M}$ if $U_{M}\supset \breve
B_{\rho}(x_{0})$ for some $\rho>0$ and if for every sequence $\rho_{j}\downarrow 0$ there is a subsequence
$\rho_{j^{\prime}}$ such that $\eta_{x_{0},\rho_{j}}M\to \C$ in $ \R^{N}\setminus \{0\}$ (convergence in the
sense of~1.4) for some cone $\C\in{\cal{}C}$ (${\cal{}C}$ denoting the set of cones in ${\cal{}P}$ as discussed
above), and also $\rho_{j^{\prime}}^{2}q(x_{0}+\rho_{j^{\prime}}x)\to r^{-2}q_{\Sigma}(x)$ ($\Sigma=\C\cap
S^{N-1}$), uniformly on compact subsets of~$\C$ (in the sense described in the remark following~1.4), where
$q_{\Sigma}$ is a non-negative locally bounded measurable function on~$\Sigma$, $\limsup_{j^{\prime}\to
\infty}\rho_{j^{\prime}}^{-n/(n+\alpha)}\||a|^{1/2}+|b|\|_{L^{n+\alpha}(\breve
B_{\rho_{j^{\prime}}}(x_{0})\setminus B_{\rho_{j^{\prime}}/2}(x_{0}))}<\infty$ and
$\rho_{j^{\prime}}^{2}|a(x_{0}+\rho_{j^{\prime}}x)|+\rho_{j^{\prime}}|b(x_{0}+ \rho_{j^{\prime}}x)|\to 0$
uniformly on compact subsets of $\C$.
 
(b) Similarly $(M, q,a,b)$ is asymptotically conic at $\infty$ if $U_{\! M}\supset \R^{N}\setminus B_{R_{0}}$ for
some $R_{0}>0$ and if for every sequence $R_{j}\uparrow \infty$ there is a subsequence $R_{j^{\prime}}$ such
that $\eta_{0,R_{j^{\prime}}}M\to \C$ in $ \R^{N}\setminus \{0\}$ (convergence in the sense of~1.4) for some
cone $\C\in{\cal{}C}$, and also $R_{j^{\prime}}^{2}q(R_{j^{\prime}}x)\to r^{-2}q_{\Sigma}(x)$ (again
$\Sigma=\C\cap S^{N-1}$), uniformly on compact subsets of~$\C$ (again in the sense described in the remark
following~1.4), where $q_{\Sigma}$ is a non-negative locally bounded measurable function on~$\Sigma$,
$\limsup_{j^{\prime}\to \infty}R_{j^{\prime}}^{-n/(n+\alpha)}\||a|^{1/2}+|b|\|_{L^{n+\alpha}(M\cap
B_{R_{j^{\prime}}}\setminus B_{R_{j^{\prime}}/2})}<\infty$ and
$R_{j^{\prime}}^{2}|a(R_{j^{\prime}}x)|+R_{j^{\prime}}|b(R_{j^{\prime}}x)|\to 0$ uniformly on compact subsets of
$\C$.

\smallskip

Notice that the definition here allows the possibility that the cone $\C$ may not be unique; that is, we may get
different cones by taking different sequences $\rho_{j},\,\rho_{j^{\prime}}$ in case~(a) and different sequences
$R_{j},\,R_{j^{\prime}}$ in case~(b). Any such cone $\C$ is called a tangent cone of $M$ (``tangent cone at
$x_{0}$'' in case~(a) and ``tangent cone at $\infty$'' in case~(b)).

We let ${\cal{}C}(M,x_{0})$ denote the (compact) set of all cones $\C\in{\cal{}C}$ which arise as in 1.12(a),(b)
according as $x_{0}\in \sing M$ or $x_{0}=\infty$ respectively,  set 
$$%
{\cal{}C}_{0}={\cal{}C}(M,x_{0}), \quad {\cal{}E}_{0}=\{\Sigma=\C\cap S^{N-1}:\C\in {\cal{}C}_{0}\}, %
$$%
and 
 $$%
 \lambda_{1}(M,x_{0})=\lambda_{1}({\cal{}E}_{0})  %
\leqno{1.13}
 $$%
 with $\lambda_{1}({\cal{}E}_{0})$ as in~1.8 with ${\cal{}C}_{0}={\cal{}C}(M,x_{0})$.  Then we have a following
corollary of Theorem~1 which, in view of the definition~1.12 of asymptotically conic, follows directly by
applying Theorem~1 iteratively in the case when ${\cal{}C}_{0}={\cal{}C}(M,x_{0})$ and when the $q_{\Sigma}$
corresponding to $\C\in{\cal{}C}_{0}$ are the functions obtained as in 1.12.

 \begin{state}{\bf{}Theorem 2.} %
   Suppose 1.1--1.4 hold and $p \in [1, \,{n\over n-2})$.  If $M\in {\cal{}P}$, $a,q:M\to \R$ with $q\ge 0$ a.e.,
$b:M\to \R^{N}$, and either $x_{0}\in \sing M$ or $x_{0}=\infty$, and $(M,q,a,b)$ is asymptotically conic at
$x_{0}$ in the sense of~1.12(a) in case $x_{0}\in\sing{}M$ and in the sense of~1.12(b) in case $x_{0}=\infty$, if
$\gamma<\gamma_{0} \equiv {n-2\over 2} -\sqrt{({n-2\over 2})^{2}+ \lambda_{1}(M,x_{0})}$, where
$\lambda_{1}(M,x_{0})$ is as in~1.13, and if $u\in W^{1,2}_{\text{\emph{loc}}}(M)\setminus\{0\}$ is a
non-negative supersolution of the equation $\Delta_{\!  M}u+b\cdot \nabla_{\! M} u +(q +a )u= 0$ on $M$, then
there is $r_{0}=r_{0}(p,\gamma,M,q, a,b)>0$ such that
\begin{align*}%
  \|u\|_{L^{p}(M\cap B_{r}\setminus B_{r/2})} %
  & \le r^{-\gamma} \text{ for all $r\ge r_{0}$ in case $x_{0}=\infty$} \\ %
  \|u\|_{L^{p}(M\cap B_{r}(x_{0})\setminus B_{r/2}(x_{0}))} %
  & \ge r^{-\gamma} \text{ for all $r\le r_{0}$ in case $x_{0}\in\sing M$};
\end{align*}%
here we continue to use the notational convention that
$\|f\,\|_{L^{p}(\Omega)}=\bigl(|\Omega|^{-1}\int_{\Omega}|f\,|^{p}\bigr)^{1/p}$, with
$|\Omega|={\cal{}H}^{n}(\Omega)$.
\end{state}
 
{\bf{}Remark:} We emphasize again that there are no a-priori continuity or indeed integrability assumptions on
$u$; $u$ is merely assumed to be non-negative a.e.\ and in $W^{1,2}_{\text{loc}}(M)\setminus\{0\}$ on the open
manifold $M$ and to be a supersolution of the equation $\Delta_{\!  M}u+b\cdot \nabla_{\! M} u +(q +a )u= 0$
locally weakly in $M$. Of course, as the above discussion and the statement of the theorem already indicates,
we can prove that $u$ automatically has integrability properties (for example part of the conclusion of the
theorem is that $u$ is automatically in $L^{p}(M\cap B_{r}\setminus B_{r/2})$ if $p<\frac{n}{n-2}$).

\smallskip

{\bf{}Proof of Theorem~2:} In case $x_{0}\in\sing M$, using the definition 1.12(a) we see that the hypotheses of
Theorem~1 are satisfied, with $\eta_{x_{0},\rho}M$ in place of $M$ and $u\circ \eta^{-1}_{x_{0},\rho}$ in place
of $u$, for all $\rho\le r_{0}$, where $r_{0}=r_{0}(p,\gamma,M,q, a,b)>0$ and where
${\cal{}C}_{0}={\cal{}C}(M,x_{0})$, the set of tangent cones of $M$ at $x_{0}$ as in the discussion preceding the
statement of Theorem~2. Thus, by Theorem~1, $\|u\|_{L^{p}(M\cap B_{\tau\rho}\setminus B_{\tau\rho/2})}
\ge \tau^{-\gamma} \|u\|_{L^{p}(M\cap B_{\rho}\setminus B_{\rho/2})}$, and, taking the choice
$\rho=r_{0}\tau^{j}$ we obtain
$$%
\|u\|_{L^{p}(M\cap B_{\tau^{j+1}r_{0}}\setminus B_{\tau^{j+1}r_{0}/2})} \ge \tau^{-\gamma}
\|u\|_{L^{p}(M\cap B_{\tau^{j}r_{0}}\setminus B_{\tau^{j}r_{0}/2})},\quad j=0,1,2,\ldots,
$$%
so by iteration we obtain the asymptotic stated in the theorem. The proof in case $x_{0}=\infty$ is a similar
iterative application of Theorem~1.

\medskip

Notice in particular that the above theorem with $p=1$ implies:

\begin{state}{\bf{}Theorem 3.} % 
  Suppose $M\in {\cal{}P}$, $x_{0}\in \sing M$, $M,a,b,q$ is asymptotically conic at $x_{0}$ in the sense
of~1.12(a) and suppose there exists $\rho>0$ such that $u$ is a non-negative supersolution of the equation
$\Delta_{\!  M}u+b\cdot \nabla_{\! M} u +(q +a )u= 0$ in $M\cap \breve B_{\rho}(x_{0})$ with $\sup_{\tau\le
\rho}\tau^{-n}\int_{M\cap B_{\tau}(x_{0})\setminus B_{\tau/2}(x_{0})}u<\infty$.  Then $\liminf_{r\downarrow
0}r^{2-n}\int_{M_{r,\delta}}q=0$ for each $\delta>0$, where $M_{r,\delta} =\{x\in M\cap B_{r}(x_{0})\setminus
B_{r/2}(x_{0}):\dist(x,\sing M)\ge \delta r\}$.
\end{state}

{\bf{}Remark:} Thus in particular there cannot exist a bounded non-negative $W^{1,2}_{\text{loc}}(M\cap
B_{\rho}(x_{0}))$ supersolution of the equation $\Delta_{\! M}u+b\cdot \nabla_{\! M} u +(q +a )u= 0$ if the
function $q$ satisfies $\liminf_{r\downarrow 0}r^{2-n}\|q\|_{L^{1}(M_{r,\delta})}>0$ for some $\delta>0$.
 
\section{Some preliminaries concerning the class \bm{${\cal{}P}$}}

First we claim there are constants $\beta_{1}= \beta_{1}({\cal{}P}),\,\beta_{2}=\beta_{2}({\cal{}P},\theta) >0$
such that
$$%
\beta_{1}\rho^{n}\le {\cal{}H}^{n}(P\cap B_{\rho}(y))
\text{ and }{\cal{}H}^{n}(P\cap B_{\theta\rho}(y))\le \beta_{2}\rho^{n}
 \leqno{2.1} %
$$%
for each $P\in {\cal{}P}$, $y\in \overline P$, $\theta\in (0,1)$ and $\rho>0$ with $B_{\rho}(y)\subset U_{P}$. 
The right inequality is in fact a direct consequence of the scale-invariance~1.2 and the first property in~1.4, and
in view of~1.3 we then have that if $P_{k}$ is an arbitrary sequence in ${\cal{}P}$ with $U_{P_{k}}\supset
B_{1}$ for each $k$ and if $P_{k}\to P\in{\cal{}P}$ in $\breve B_{1}$ with $U_{P}\supset \breve B_{1}$ then
${\cal{}H}^{n}(\breve B_{\theta}\cap \{x\in P:\dist(x,\sing P)<\delta\})\le C\delta^{2}$ for each $\theta\in [1/2,1)$. 
In view of the $C^{1}$ and Hausdorff distance sense convergence of~1.4, it evidently follows that
${\cal{}H}^{n}\res P_{k}\to {\cal{}H}^{n}\res P \textrm{ in }\breve B_{1}$; that is $\int_{P_{k}}f\,d{\cal{}H}^{n}\to
\int_{P}f\, d{\cal{}H}^{n}$ for each fixed continuous $f$ with compact support in $\breve B_{1}$. Thus we have
established
\begin{align*}%
\tag*{2.2}  %
&P,P_{k}\in{\cal{}P} \text{ with } U_{P_{k}}\supset \breve B_{1} \text{ and } P_{k}\to P
                                                                   \text{ in the sense of~1.4 in  }\breve B_{1} \\
&\hskip2.5in     \Rightarrow  {\cal{}H}^{n}\res P_{k}\to {\cal{}H}^{n}\res P \text{ in }\breve B_{1}. %
\end{align*}%
To prove the left inequality in~2.1, suppose on the contrary that $\rho_{k}^{-n}{\cal{}H}^{n}(P_{k}\cap{}
B_{\rho_{k}}(y_{k}))<k^{-1}$, $k=1,2,\ldots$, with $y_{k}\in\overline P_{k}$ and $U_{P_{k}}\supset
B_{\rho_{k}}(y_{k})$. Then, with $\wtilde P_{k}=\eta_{y_{k},\rho_{k}}P_{k}$, we have by~1.2 that $\wtilde
P_{k}\in {\cal{}P}$ with $U_{\wtilde P_{k}}\supset B_{1}$, $0\in$ the closure of $\wtilde P_{k}$, and
${\cal{}H}^{n}(\wtilde P_{k}\cap B_{1})<k^{-1}$ for each $k$. Then by~1.4 there is $\wtilde P\in{\cal{}P}$ and a
subsequence $\wtilde P_{k^{\prime}}\to \wtilde P$ in $\breve B_{1}$ with $U_{\wtilde P}\supset \breve B_{1}$,
with ${\cal{}H}^{n}(\wtilde P)=0$ (by~2.2) and with $0$ in the closure of $\wtilde P$ (by the Hausdorff distance
convergence), contradicting the assumption that all elements of ${\cal{}P}$ are $n$-dimensional
submanifolds.

Notice that if $\theta\in [{1\over 2}, 1)$ is given, we can now bound the number of connected components
$P_{\ast}$ of $P\cap \breve B_{\rho}(y)$ which intersect $B_{\theta\rho}(y)$ in case $U_{P}\supset B_{\rho}(y)$. 
Indeed, since for each such component $P_{\ast}$ we have $z\in{}P_{\ast}\cap{}B_{\theta\rho}(y)$, and hence
$\breve B_{{1\over{}2}(1-\theta)\rho}(z)\subset\breve B_{\wtilde\theta\rho}(y)\subset \breve B_{\rho}(y)\subset
U_{P_{\ast}}$, where $\wtilde\theta={1+\theta\over{}2}$, the left inequality in~2.1 gives
 $$%
 {\cal{}H}^{n}(P_{\ast}\cap B_{\wtilde\theta\rho}(y)) \ge \beta_{1}\rho^{n},  %
\leqno{2.3}
 $$%
 for suitable $\beta_{1}=\beta_{1}(\theta,{\cal{}P})$ whereas the sum of ${\cal{}H}^{n}(P_{\ast}\cap
B_{\wtilde\theta\rho}(y))$ over all such components $P_{\ast}$ is $\le {\cal{}H}^{n}(P\cap
B_{\wtilde\theta\rho}(y))\le \beta_{2}\rho^{n}$ for some $\beta_{2}=\beta_{2}(\theta,{\cal{}P})$ by the right
inequality in~2.1, whence the number $Q$ of such components satisfies
 $$%
 Q\le 1 + \beta_{1}^{-1}\beta_{2}. %
 \leqno{2.4} %
 $$%
 Finally, we show that the conditions 1.1--1.4 are sufficient to give a ``restricted Poincar\'e type'' inequality
on each $P\in {\cal{}P}$:

\begin{state}{\bf 2.5 Theorem.}  %
Let ${\cal{}P}$ satisfy the conditions 1.1--1.4.  Then for each $\theta_{0}\in
[\ha, 1)$ there are constants $C= C({\cal{}P},\, \theta_{0})>0$, $\delta=\delta({\cal{}P},\theta_{0})\in(0, \ha]$
such that
$$%
\Bigl(\int_{P\cap B_{\theta_{0}}}\varphi^{\kappa}\Bigr)^{1/\kappa} \le C\int_{P\cap \breve B_{1}}|\nabla
\varphi|,\qquad \kappa = {n\over n-1},
$$%
whenever $P\in {\cal{}P}$ with $U_{P}\supset \breve B_{1}$ and $\varphi$ is a non-negative $C^{1}(P\cap
\breve B_{1})$ function satisfying the inequality ${\cal{}H}^{n}(\support \varphi)<\delta$.
\end{state}

{\bf 2.6 Remarks: (1)} An examination of the proof will show that for this lemma it would suffice that
${\cal{}H}^{n-1}(\sing P)=0$ for each $P\in {\cal{}P}$ in place of the condition~1.3.

{\bf{}(2)} By replacing $\varphi$ by $|\varphi|^{2(n-1)/(n-2)}$ for $n\ge 3$ and by $\varphi^{2q}$ for arbitrary
$q>1$ in case $n=2$, and using the H\"older inequality on the right side, we see that the inequality of~2.5 admits
the squared version
$$%
\Bigl(\int_{P\cap B_{\theta_{0}}}\varphi^{2\kappa}\Bigr)^{1/\kappa} \le  %
                                                       C\int_{P\cap \breve B_{1}}|\nabla \varphi|^{2}, %
$$%
with $\kappa=n/(n-2)$ and $C=C({\cal{}P},N)$ in case $n\ge 3$, and in case $n=2$ the same with arbitrary
$\kappa>1$. (Of course we still require the restriction ${\cal{}H}^{n}(\support \varphi)<\delta$ here.)

Before we begin the proof of~2.5 we observe that, using it in combination with a partition of unity for
$B_{1}$ consisting of smooth functions, each of which has support in a set of diameter $\le \delta$, we
conclude the following.

\begin{state}{\bf 2.5$^{\prime}$  Corollary.} %
 If the hypotheses are as in~2.5, except that we drop the condition
that ${\cal{}H}^{m}(\support \varphi)<\delta$, then
 $$%
 \Bigl(\int_{P\cap B_{\theta_{0}}}\varphi^{\kappa}\Bigr)^{1/\kappa} \le  %
                         C\int_{P\cap \breve B_{1}}\Bigl(|\nabla \varphi|+|\varphi|\Bigr). %
 $$%
\end{state}

{\bf{}2.6$^{\prime}$ Remark:} As in Remark 2.6(2), there is a squared version of the above inequality:
$$%
\Bigl(\int_{P\cap B_{\theta_{0}}}\varphi^{2\kappa}\Bigr)^{1/\kappa} \le  %
                               C\int_{P\cap \breve B_{1}}(|\nabla \varphi|^{2}+|\varphi|^{2}), %
$$%
with $\kappa=n/(n-2)$ and $C=C({\cal{}P},N)$ in case $n\ge 3$, and in case $n=2$ the same with arbitrary
$\kappa>1$.

{\bf Proof of Theorem~2.5}: It is a well-known consequence of the coarea formula and the fact that
${\cal{}H}^{n-1}(\sing P\cap \breve B_{1})=0$ that such an inequality is equivalent to the fact that $\exists \,
C=C({\cal{}P},\theta)>0$ such that
$$% 
{\cal{}H}^{n}(Q\cap B_{\theta_{0}})^{1/\kappa} \le C {\cal{}H}^{n-1}(\partial Q\cap \breve B_{1})
$$%  
whenever $Q$ is a relatively open subset of $P\cap \breve B_{1}$ with ${\cal{}H}^{n}(Q) <\delta$, with
boundary $\partial Q=\overline Q\setminus Q\, (\subset\overline P)$ such that $\partial Q\cap P\cap \breve
B_{1}$ is locally $C^{1}$.

Take $\delta\in (0,{1\over{}4}]$ such that $\delta^{1/2}$ is smaller than the volume $\omega_{n}$ of the unit ball
in $\R^{n}$ and also smaller than the constant $\beta_{1}$ in~2.1, and assume (to get a contradiction) that
$P_{k}$ is a sequence in ${\cal{}P}$ such that for each $k$ there is a relatively open subset $Q_{k}\subset
P_{k}\cap B_{1}$ with ${\cal{}H}^{n}(Q_{k})<\delta$ and with boundary $\partial Q_{k}\cap \breve B_{1}$ such
that ${\cal{}H}^{n-1}(\partial Q_{k}\cap \breve B_{1}) < {1\over{}k}({\cal{}H}^{n}(Q_{k}\cap
B_{\theta_{0}}))^{1/\kappa}\to 0$.  Now for each point $y$ in $Q_{k}\cap B_{\theta_{0}}$,
$\rho^{-n}{\cal{}H}^{n}(Q_{k}\cap B_{\rho}(y))$ has limiting value~$\omega_{n}>\delta^{1/2}$ as
$\rho\downarrow 0$ and has value $\le C(\theta_{0})\delta<\delta^{1/2}$ when $\rho={1\over{}2}(1-\theta_{0})$,
assuming $\delta$ small enough (depending on $\theta_{0}$), and
$\lim_{\sigma\downarrow\rho}\sigma^{-n}{\cal{}H}^{n}(Q_{k}\cap
B_{\sigma}(y))=\rho^{-n}{\cal{}H}^{n}(Q_{k}\cap B_{\rho}(y))$ for each $\rho\in (0,{1\over{}2}(1-\theta_{0}))$, so
there is a smallest value $\rho(y,k)$ of $\rho\in (0,{1\over{}2}(1-\theta_{0}))$ such that ${\cal{}H}^{n}(Q_{k}\cap
B_{\rho}(y))\le \delta^{1/2} \rho^{n}$.  Thus
$$%
{\cal{}H}^{n}(Q_{k}\cap B_{\rho_{y,k}}(y))= \delta^{1/2} \rho_{y,k}^{n} \text{ and } %
{\cal{}H}^{n}(Q_{k}\cap B_{\rho}(y))\ge \delta^{1/2}\rho^{n}\,\,\forall\, \rho\in (0,\rho_{y,k}]  %
\leqno{(1)}
$$%
By the Besicovich covering lemma there is a subcollection $\{B_{\rho_{j,k})}(y_{j,k})\}$ of such balls which
covers $Q_{k}$ and which decomposes into a fixed number $J=J(N)$ of pairwise-disjoint subcollections.  For
each $k$ we must then have at least one of these balls, say $B_{\rho_{k}}(y_{k})$, with
$$%
{\cal{}H}^{n-1}(B_{\rho_{k}}(y_{k})\cap\partial{}Q_{k}\cap \breve B_{1})  %
                   \le k^{-1/2}({\cal{}H}^{n}(Q_{k}\cap{}B_{\rho_{k}}(y_{k})))^{1/\kappa}, %
\leqno{(2)}
$$%
because otherwise we would have 
$$%
   k^{-1/2}({\cal{}H}^{n}(Q_{k}\cap B_{\rho_{j,k}}(y_{j,k})))^{1/ \kappa}   %
                                     <{\cal{}H}^{n-1}(B_{\rho_{j,k}}(y_{j,k})\cap \partial Q_{k}\cap \breve B_{1}) %
$$%
for each $j$, and we could sum over $j$ to conclude that
\begin{align*}%
\bigl({\cal{}H}^{n}(Q_{k}\cap B_{\theta_{0}})\bigr)^{1/\kappa} %
&\le \bigl(\sum_{j}{\cal{}H}^{n}(Q_{k}\cap B_{\rho_{j,k}}(y_{j,k}))\bigr)^{1/ \kappa}\le   %
                 \sum_{j}\bigl({\cal{}H}^{n}(Q_{k}\cap B_{\rho_{j,k}}(y_{j,k}))\bigr)^{1/ \kappa}  \\  %
&\le k^{1/2}\sum_{j}{\cal{}H}^{n-1}(B_{\rho_{j,k}}(y_{j,k})\cap \partial Q_{k}\cap \breve B_{1}) %
                                                             \le k^{1/2}J {\cal{}H}^{n-1}(\partial Q_{k}\cap \breve B_{1}) %
\end{align*}%
contrary to the original choice of $Q_{k}$ for sufficiently large $k$. (Notice that here we use the inequality
$(\sum_{j} a_{j})^{1/\kappa}\le \sum_{j}a_{j}^{1/\kappa}$.)

Now with $y_{k},\rho_{k}$ as in~(2), let $Q^{\prime}_{k}\equiv \eta_{y_{k}, \rho_{k}}Q_{k}$, and
$P^{\prime}_{k}\equiv \eta_{y_{k}, \rho_{k}}P_{k}$. By the compactness~1.4 we have a subsequence of
$P^{\prime}_{k}\to P$ in $\breve B_{1}$, where $P\in{\cal{}P}$ with $U_{P}\supset\breve{B}_{1}$.  Let
$\zeta=(\zeta^{1},\ldots,\zeta^{N})$ be a fixed $C_{\text{c}}^{\infty}(\breve B_{1};\R^{N})$ function with support
$\zeta|P$ contained in a compact subset $K$ of $P\cap \breve B_{1}$, and let $\widetilde Q_{k}=
\Psi_{k}(Q_{k}^{\prime})$, as in the remark following~1.4, be $C^{1}$ on an open set $U$ containing $K$ and
satisfy $\wtilde Q_{k}\equiv \Psi_{k}(Q_{k}^{\prime})\subset P$ and $\|\Psi_{k}-I\|_{C^{1}(U)}\to 0$ as $k\to
\infty$.  Then ${\cal{}H}^{n-1}(\partial\widetilde{Q}_{k}\cap \breve B_{1})\to{}0$ and hence
$\int_{\widetilde{Q}_{k}}\dvg_{P}\zeta\to 0$, so that by the BV compactness theorem and the arbitrariness of
$B_{\tau}(y)$ there is a measurable $Q\subset P$ and a subsequence of $\wtilde Q_{k}$ such that the indicator
functions $\chi_{\wtilde Q_{k}}$ converge strongly in $L^{1}$ on $P\cap \breve B_{1}$ to $\chi_{Q}$; then
$\int_{Q}\dvg_{P} \zeta =\lim\int_{\wtilde Q_{k}} \dvg_{P}\zeta=0$, which means that the indicator function of
$Q$ is locally constant in $P\cap \breve B_{1}$. Thus, up to a set of measure zero, $Q$ is a union of
components of $P\cap \breve B_{1}$. By virtue of~2.4 we have that at most finitely many components of $P\cap
\breve B_{1}$ can intersect the ball $B_{1/2}$, and hence at most finitely many of the components of $P\cap
\breve B_{1}$ which comprise $Q$ can intersect $B_{1/2}$.  On the other hand by construction we arranged
that $\tau^{-n} {\cal{}H}^{m}(Q\cap B_{\tau})\in [\delta^{1/2},\infty)$ for each $\tau<1$, and hence there is at
least one component $\widetilde P$ of these finitely many components containing $0$ in its closure. That is,
there is a component $\widetilde P$ of $P\cap \breve B_{1}$ with $\widetilde P\subset Q\cap \breve B_{1}$ (up
to a set of measure zero) and $0\in$ the closure of $\widetilde P$.  But also by construction we have
${\cal{}H}^{m}(Q\cap \breve B_{1}) \le \delta^{1/2}$, which means that ${\cal{}H}^{m}(\widetilde P)\le
\delta^{1/2}$.  Since $\widetilde P\in {\cal{}P}$ with $U_{\widetilde P}=\breve B_{1}$ (by the reducibility
hypothesis~1.1) with $0\in \text{ closure }\widetilde P$, this contradicts the bounds~2.1 since
$\delta^{1/2}<\beta_{1}$, where $\beta_{1}$ is as in~2.1. Thus the proof is complete.

\section{Stability Inequality}

Here $M$ will denote a fixed element of $\mathcal{P}$ and $q$ will be a non-negative locally bounded
measurable function on $M$.  $u$ will denote a positive $W^{1,2}_{\text{loc}}(M)$ supersolution of the
equation $\Delta_{M}u+ q\, u=0$ with $q$ non-negative measurable and locally bounded on $M$. This means
that
$$%
\int_{M}(-\nabla u\cdot \nabla \zeta + q\,\zeta u)\le 0,\quad \zeta\in C^{1}_{c}(M)\text{ with }\zeta\ge 0. 
\leqno{3.1}
$$% 
(Thus in the present section there are no perturbation terms $a,b$ as in~1.9.)  We claim that then we have the
``stability inequality''
$$%
\int_{M}(|\nabla \zeta|^{2}-q\,\zeta^{2})\ge 0,\quad \zeta\in C_{c}^{1}(U_{M}), 
\leqno{3.2}
$$%
and in particular that $\int_{M\cap K} q<\infty$ for each compact $K\subset U_{M}$.  Notice that while the
definition~3.1 requires $\zeta$ to vanish in a neighborhood of the singular set, the inequality~3.2 does not.  To
prove~3.2, first take any non-negative $\zeta\in C^{1}_{c}(U_{M})$ and let $s_{\delta} : M\to [0,\,1]$ be a
smooth function with compact support in $\support\zeta\cap M$ defined as follows: First use the definition of
finite $\mathcal{H}^{n-2}$-measure and the compactness of $\support \zeta\cap \sing M$ there is a constant
$\beta$ such that for each $\delta\in (0,1)$ we can select a finite cover $\breve B_{\rho_{j}/2}(y_{j}),j=1,\ldots,
Q$ of $\text{sing}M\cap \support \zeta$ by balls with centers in $\sing M$, $\sup_{j}\rho_{j}<\delta$ and
$\sum_{j}\rho_{j}^{n-2}<\beta$.  Next, for $j=1,\ldots, Q$, let $\zeta_{j}$ be a smooth function on $M$ with
$\zeta_{j}\equiv 0$ on $B_{\rho_{j}/2}(y_{j})$, $\zeta_{j}\equiv 1$ on $B_{\rho_{j}}(y_{j})$, $0\le \zeta_{j}\le 1$
everywhere, and $|\nabla \zeta_{j}|\le 3\rho_{j}^{-1}$.  Then with $s_{\delta}=\min\{\zeta_{1},\ldots, \zeta_{Q}\}$
we have $s_{\delta}(x)\equiv 1$ for $\text{dist}(x,\sing M)>\delta$, $s_{\delta}\equiv 0$ in a neighborhood of
$\sing M\cap B_{\rho}(y)$, while $\int_{M\cap \support\zeta}|\nabla s_{\delta}|^{2}\le C \sum_{j}\rho_{j}^{n-2} \le
C\beta $, where $C\le \beta_{2}$ with $\beta_{2}$ as in~2.1. Now use~3.1 with $(\epsilon+u)^{-1}\zeta^{2}\,
s^{2}_{\delta}$ in place of $\zeta$. Since $\zeta\,s_{\delta}$ has compact support in $M$ this is a legitimate
choice, and~3.1 gives
$$%
\int_{M}\zeta^{2}\, s^{2}_{\delta}({u\over u+\epsilon}q +|\nabla w|^{2})\le   %
                                   -2\int_{M}(s_{\delta}\zeta\nabla w\cdot \nabla (\zeta s_{\delta}))  %
$$%
with $w=\log (\epsilon+u)$.  Using the Cauchy inequality $a\cdot b\le |a|^{2}+{1\over 4}|b|^{2}$ on the right side,
we thus deduce that $\int_{M}\zeta^{2}\, s^{2}_{\delta}{u\over u+\epsilon}q +|\nabla w|^{2}\le \int_{M}|\nabla
(s_{d}\zeta)|^{2}\le C$ with constant $C$ independent of $\delta$, so that by letting $\delta\downarrow 0$ we
have $\int_{M}\zeta^{2}|\nabla w|^{2}<\infty$. On the other hand~3.1 implies
$$%
\int_{M}\zeta^{2}\, s^{2}_{\delta}({u\over u+\epsilon}q +|\nabla w|^{2})\le   %
                          -2\int_{M}(s_{\delta}^{2}\zeta\nabla w\cdot \nabla \zeta )  %
                                         -2\int_{M}(s_{\delta}\zeta^{2}\nabla w\cdot \nabla s_{\delta} ). %
$$%
Letting $\delta\downarrow 0$ and using Cauchy-Schwarz to check that the last integral on the right tends to
zero, and we then have
$$%
\int_{M}\zeta^{2}\, ({u\over u+\epsilon}q +|\nabla w|^{2})\le   %
                          -\int_{M}((2\zeta\nabla w)\cdot \nabla \zeta )  %
$$%
and so, letting $\epsilon\downarrow 0$ and using $ab\le \tfrac{1}{4} a^{2}+b^{2}$, we conclude~3.2 as claimed.

\section{Compact classes of cones}

The present section will be needed in the proof of Theorem~2 because we do not assume in the definition~1.12
of asymptotically conic that $M\in {\cal{}P}$ necessarily has a unique tangent cone at points $x_{0}\in\sing M$
or at $\infty$. To overcome this difficulty we shall use that fact if $M$ is as in~Theorem~2 then the set of all
possible tangent cones $\C$ of $M$ arising as in~1.12 (at a point $x_{0}\in \sing M$ in case 1.12(a) or at $\infty$
in case 1.12(b)) is a compact subfamily of ${\cal{}P}$ with respect to the natural Hausdorff distance metric
$d_{1}$ defined below in~4.1.

Let ${\cal{}C}$ denote the set of all cones in~${\cal{}P}$ as in~\S1. As we mentioned in~\S1, the Hausdorff
distance metric $d$ on ${\cal{}E}=\{\Sigma=\C\cap S^{N-1}:\C\in {\cal{}C}\}$ makes ${\cal{}E}$ into a compact
metric space, and of course we can metrize ${\cal{}C}$ by the metric $d_{1}$ given by
$d_{1}(\C_{1},\C_{2})=d(\Sigma_{1},\Sigma_{2})$, where $\Sigma_{j}=\C_{j}\cap S^{N-1}$, and then
$$% 
\mbox{ ${\cal{}C}$, equipped with the metric $d_{1}$, is a compact metric space.} %
\leqno{4.1} %
$$% 
Now as in~\S1 let ${\cal{}C}_{0}$ be any fixed compact subset of ${\cal{}C}$, let ${\cal{}E}_{0}=\{\C\cap
S^{N-1}:\C\in{\cal{}C}_{0}\}$, and for each $\Sigma\in{\cal{}E}_{0}$ assume we have a non-negative
$q_{\Sigma}$ such that the collection ${\cal{}Q}_{0}$ of all such $q_{\Sigma}$ satisfies the compactness
assumption of~1.6.  In view of the compactness of ${\cal{}C}_{0}$, we must then have for each $\tau >0$
$$% 
q_{\Sigma}(\omega)\le \Lambda_{{\cal{}C}_{0},\tau},\quad \omega\in \Sigma_{\tau}, %
\leqno{4.2}%
$$%
where $\Sigma_{\tau} = \{x\in \Sigma: \dist(x,\sing \Sigma)>\tau\}$, and $\Lambda_{{\cal{}C}_{0},\tau}$ is
a fixed constant depending only on ${\cal{}C}_{0}$ and $\tau$, and not depending on the particular cone $\C$.

For $\Sigma\in{\cal{}E}_{0}$ we continue to define $\lambda_{1}(\Sigma)$ as in~1.7.  We observe that for each
$\tau\in (0,{1\over{}2}]$ we have $\epsilon(\tau)=\epsilon(\tau,\Sigma,q_{\Sigma})\downarrow 0$
as $\tau \downarrow0$ such that
\begin{align*}%
  \tag*{4.3} %
  &  \text{  if } \Sigma\in {\cal{}E}_{0} \text{ is such that $\exists$ a non-negative $u\in
W^{1,2}_{\text{loc}}(\C\cap \breve B_{1}\setminus B_{\tau})\setminus\{0\}$ with } \\   
\noalign{\vskip-4pt} %
  &\hskip.2in \text{ $\Delta_{\C}u+r^{-2} q_{\Sigma}u\le 0$ weakly on $\C\cap \breve B_{1}\setminus
B_{\tau}$, then } %   
\lambda_{1}(\Sigma) \ge -\textstyle{\Bigl({n-2\over 2}\Bigr)^{2}}-\epsilon(\tau). %
\end{align*}%
To prove this we use the stability inequality~3.2 with $\zeta(x) = \zeta_{1}(r)\zeta_{2}(\omega)$, where $r=|x|$,
$\omega= |x|^{-1}x$, $\zeta_{1}\in C_{c}^{1}(\tau,1)$, and $\zeta_{2}\in C_{c}^{1}(\R^{N}\setminus \{0\})$
homogeneous of degree zero with $\support\zeta_{2}\cap \sing\Sigma=\emptyset$.  Then~3.2 implies
$$%
0\le \int_{0}^{1}(\zeta_{1}^{\prime}(r))^{2}\, r^{n-1} dr
\int_{\Sigma}\zeta_{2}^{2}(\omega)\,d\omega + \int_{0}^{1}\zeta_{1}^{2}(r)\, r^{n-3}dr \int_{\Sigma}(|\nabla
\zeta_{2}|^{2}- q_{\Sigma} \zeta_{2}^{2})\,d\omega,
$$%
whence, taking inf over all $\zeta_{1},\, \zeta_{2}$ with $L^{2}$ norms equal to~1, we conclude
that
$$%
\lambda_{1}(\Sigma) \ge - \inf {\int_{0}^{1}(\zeta_{1}^{\prime})^{2}\, r^{n-1}dr\over
\int_{0}^{1}\zeta_{1}^{2}\, r^{n-3}dr},
$$%
where the inf is over all $\zeta_{1}\in C^{1}(\tau,1)$ with compact non-empty support. It is a standard calculus
fact that if $\tau=0$ then the quantity on the right is exactly $-\left({n-2\over 2}\right)^{2}$ and hence in general
it is $\ge -\left({n-2\over 2}\right)^{2}-\epsilon(\tau)$ with $\epsilon(\tau)\downarrow 0$ as $\tau\downarrow0$. 
this gives the required inequality $\lambda_{1}(\Sigma)\ge -\bigl({n-2\over{}2}\bigr)^{2}-\epsilon(\tau)$ as
claimed.

\medskip

The following lemma ensures we can always select a collection of eigenfunctions with good positivity
properties on domains in $\Sigma\in{\cal{}E}_{0}$ and with eigenvalues not much bigger than the value
$\lambda_{1}({\cal{}E}_{0})=\sup\{\lambda_{1}(\Sigma): \Sigma\in{\cal{}E}_{0}\}$ of~1.8.
In this lemma we use the notation that
$$%
{\cal{}E}_{0}(\Lambda) =\{\Sigma\in{\cal{}E}_{0}:\lambda_{1}(\Sigma)\ge \Lambda\}
$$%
for $\Lambda\in \R$. Observe that ${\cal{}E}_{0}(\Lambda)$ is a closed (hence compact) subset of
${\cal{}E}_{0}$, which is easily checked by using the Rayleigh quotient definition~1.7 together with~1.6 and the
local $C^{1}$ convergence described in~1.4.

\begin{state}{\bf{}4.4 Lemma.} %
  For each $\delta>0,\Lambda\in \R$, $\exists\,\,\tau=\tau(\delta,\Lambda,{\cal{}E}_{0},{\cal{}Q}_{0})>0$ such that
the following holds.  For each $\Sigma\in{\cal{}E}_{0}(\Lambda)$ there are connected open $C^{1}$ domains
$\Omega_{1},\ldots,\Omega_{Q}$, $Q\le Q_{0}=Q_{0}({\cal{}E}_{0})$, with $\overline\Omega_{j}\subset \Sigma$
and $\overline\Omega_{i}\cap \overline\Omega_{j}=\emptyset\,\forall\,i\neq j$, and corresponding non-negative
$\varphi_{j,\delta}\in W^{1,2}_{0}(\Omega_{j})\cap C^{0}(\overline\Omega_{j})$ with
$$%
\max\varphi_{j,\delta}=1\text{ and }-(\Delta_{\Sigma}\varphi_{j,\delta}+
                               q_{\Sigma}\varphi_{j,\delta})=\lambda_{1}(\Omega_{j})\varphi_{j,\delta}%  
\text{ weakly on }\Omega_{j}, %
$$%
where 
$$%
\lambda_{1}(\Omega_{j})\le \lambda_{1}({\cal{}E}_{0})+\delta \text{ for each }j=1,\ldots,Q
$$%
and 
$$%
{\cal{}H}^{n}(\Sigma\setminus(\cup_{j}\{x\in\Omega_{j}:\varphi_{j,\delta}(x)>\tau\}))<\delta.
$$%
\end{state}%

{\bf{}Remark:} An essential feature of the above lemma is that the constant $\tau$ does not depend on the
particular $\Sigma\in{\cal{}E}_{0}(\Lambda)$ under consideration, so $\tau$ is chosen and then the lemma applies
uniformly across the whole class ${\cal{}E}_{0}(\Lambda)$ and corresponding ${\cal{}Q}_{0}$ (as in~1.6).

{\bf{}Proof of Lemma~4.4:} If $\Sigma\in{\cal{}E}_{0}(\Lambda)$ with connected components
$\Sigma_{1},\ldots,\Sigma_{Q}$ (so that $Q\le Q_{0}=Q_{0}({\cal{}E}_{0})$ by~2.4), then for each sufficiently
small $\tau>0$ we can select connected open $C^{1}$ domains $\Omega_{j,\tau}$ with
$$%
\Sigma_{j,\tau} \subset \Omega_{j,\tau}\subset\overline\Omega_{j,\tau} \subset\Sigma_{j},
$$%
where we use the notation $\Sigma_{j,\tau}=\{x\in\Sigma_{j}:\dist(x,\sing \Sigma_{j})>\tau\}$.  Then we can take
$\varphi_{j}^{(\tau)}\in W_{0}^{1,2}(\Omega_{j,\tau})\setminus\{0\}$ minimizing the Rayleigh quotient
$$%
\Bigl(\int_{\Omega_{j,\tau}}\varphi^{2}\Bigr)^{-1} %
             \int_{\Omega_{j,\tau}}\bigl(|\nabla\varphi|^{2}-q_{\Sigma}\varphi^{2}\bigr) %
$$%
over $\varphi\in W_{0}^{1,2}(\Omega_{j,\tau})\setminus\{0\}$. Letting $\lambda_{j,\tau}$ denote the minimum
value, we then have, for small enough $\tau=\tau(\Sigma,\delta)$,  that $\varphi_{j}^{(\tau)}$ is
non-negative a weak solution of the equation
$$%
-(\Delta_{j}\varphi_{j}^{(\tau)}+q_{\Sigma}\varphi_{j}^{(\tau)})=\lambda_{j,\tau}\varphi_{j}^{(\tau)},
$$%
and, in the notation of~1.7, $\lambda_{j,\tau}\to \lambda_{1}(\Sigma_{j})$ as $\tau\downarrow 0$ for each
$j=1,\ldots,Q$. Also by the De~Giorgi Nash Moser theory (\cite[\S8.8--\S8.10]{GilT83}) we know that
$\varphi_{j}^{(\tau)}\in W^{1,2}_{0}(\Omega_{j,\tau})\cap C^{0}(\overline\Omega_{j,\tau})$ 
and that
$$%
\begin{aligned}%
&\sup_{\Omega_{j,\tau}}\varphi_{j}^{(\tau)}\le  %
                     C(\Sigma,\Lambda,\tau)\|\varphi_{j}^{(\tau)}\|_{L^{2}(\Omega_{j,\tau})}, \\  %
&\hskip-.55in\inf_{\{x\in\Omega_{j,\tau}:\dist(x,\partial\Omega_{j,\tau})>\tau\}}\varphi_{j}^{(\tau)} %
                                   \ge C(\Sigma,\tau)^{-1}\|\varphi_{j}^{(\tau)}\|_{L^{2}(\Omega_{j,\tau})}. %
\end{aligned}%
$$%
So we can normalize so that $\max_{\Omega_{j,\tau}}\varphi_{j}^{(\tau)} =1$ and then
$\inf_{\{x\in\Omega_{j,\tau}:\dist(x,\partial\Omega_{j,\tau})>\tau\}}\varphi_{j}^{(\tau)}\ge \sigma$ with
$\sigma=\sigma(\Sigma,\tau)>0$, hence for sufficiently small $\tau=\tau(\Sigma,q_{\Sigma},\delta)$
and small enough $\sigma=\sigma(\Sigma,\tau)>0$ we have
$$%
{\cal{}H}^{n}(\Sigma\setminus(\cup_{j=1}^{Q}\{x\in\Omega_{j,\tau}:\varphi_{j}^{(\tau)}>\tau\}))<\delta.
$$%  
Thus with such a choice of $\tau$ and $\sigma$, we use the notation
$$%
\Omega_{j}=\Omega_{j,\tau}, \quad \varphi_{j,\delta} =\varphi_{j}^{(\tau)},
$$%
and then we have the required properties stated in the lemma except that the constant $\tau$ depends on the
particular $\Sigma\in{\cal{}E}_{0}(\Lambda)$, i.e., $\tau=\tau(\Sigma,q_{\Sigma},\delta)$.  But now observe that
by~1.4 and~1.6 there is $\epsilon=\epsilon(\Sigma,q_{\Sigma})$ such that each $\wtilde\Sigma$ in the Hausdorff
distance metric ball ${\cal{}B}_{\epsilon}(\Sigma)\subset{\cal{}E}_{0}(\Lambda)$ of radius $\epsilon$ and center
$\Sigma$ is $C^{1}$ close to $\Sigma$ in a neighborhood of the compact set
$\cup_{j=1}^{Q_{0}}\overline\Omega_{j}\subset \Sigma$, and correspondingly the function $q_{\wtilde \Sigma}$
is uniformly close to $q_{\Sigma}$ in this neighborhood.  Specifically, for any $\eta>0$ and small enough
$\epsilon=\epsilon(\Sigma,q_{\Sigma},\eta,\Omega_{1},\ldots,\Omega_{Q})>0$, there is a fixed $U$, open in
$\R^{N}$, with $\cup_{j=1}^{Q}\overline\Omega_{j}\subset U$ and for each $\wtilde\Sigma\in
{\cal{}B}_{\epsilon}(\Sigma)$ there is a $C^{1}$ map $\widetilde \Psi:U\to U$ with $\|\wtilde\Psi-I\|_{C^{1}}<\eta$
and such that $\wtilde\Psi(U\cap \Sigma)=U\cap \wtilde\Sigma$, $\max_{\cup_{j=1}^{Q}}|q_{\wtilde\Sigma}\circ
\wtilde\Psi-q_{\Sigma}|<\eta$ and $\wtilde\Omega_{j}=\wtilde\Psi(\Omega_{j})$ have pairwise disjoint closures
contained in $\wtilde\Sigma$ and a of distance $<\eta$ from $\Omega_{j}$ in the Hausdorff distance sense. 
Thus, if we take $\eta=\eta(\Sigma,q_{\Sigma},\delta)$ is suitably small (technically depending also on the choice
of $\Omega_{j}$ for $\Sigma$, but those domains are determined by $\Sigma$ and $\delta$ also) then we have,
for all $\wtilde\Sigma\in {\cal{}B}_{\epsilon}(\Sigma)$, open $\wtilde\Omega_{j}\subset\wtilde \Sigma$ and
functions $\wtilde \varphi_{j,\delta}\in W_{0}^{1,2}(\wtilde \Omega_{j})\cap
C^{0}(\op{closure}\wtilde\Omega_{j})$ with
$$%
\max\wtilde\varphi_{j,\delta}=1\text{ and  }  %
     -(\Delta_{\wtilde\Sigma}\wtilde\varphi_{j,\delta}+ q_{\wtilde\Sigma}\wtilde\varphi_{j,\delta})=  %
                  \lambda_{1}(\wtilde\Omega_{j})\wtilde\varphi_{j,\delta}% 
    \text{ weakly on }\wtilde\Omega_{j}, %
$$%
where 
$$%
\lambda_{1}(\wtilde\Omega_{j})\le \lambda_{1}({\cal{}E}_{0})+2\delta \text{ for each }j=1,\ldots,Q
$$%
and 
$$%
{\cal{}H}^{n}(\wtilde\Sigma\setminus
(\cup_{j}\{x\in\wtilde\Omega_{j}:\wtilde\varphi_{j,\delta}(x)>\tau\}))<2\delta.
$$%
with constant $\tau=\tau(\Sigma,q_{\Sigma},\delta)$.  Since ${\cal{}E}_{0}(\Lambda)$ is compact, we can select
finitely many such balls ${\cal{}B}_{\epsilon(\Sigma_{k})}(\Sigma_{k}), k=1,\ldots,S$, $S=S(\Lambda,
{\cal{}E}_{0},\delta)$, such that ${\cal{}E}_{0}(\Lambda)=
\cup_{k=1}^{S}{\cal{}B}_{\epsilon(\Sigma_{k})}(\Sigma_{k})$. Taking the minimum of the corresponding
constants $\tau(\Sigma_{k},q_{\Sigma_{k}},\delta)/2, \,k=1,\ldots,S$, we thus have the conclusion stated in the
lemma.

\smallskip

{\bf{}Remarks: (1)\,} Notice that in the above proof we first selected domains
$\Omega_{j}=\Omega_{j}(\Sigma,q_{\Sigma})$ with each $\Omega_{j}$ engulfing all but a thin boundary strip of
one of the connected components of $\Sigma$, but the reader should observe that the corresponding $C^{1}$
domains $\wtilde\Omega_{j}\subset \wtilde\Sigma$ of the nearby $\wtilde\Sigma\in
{\cal{}B}_{\epsilon(\Sigma)}(\Sigma)$ may in fact be only a small fraction of one of the components of
$\wtilde\Sigma$ (e.g.\ the union of two or more of the $\wtilde\Omega_{j}$ may be needed to encompass most
of a single component of $\wtilde \Sigma$).  This is because the nearby $\wtilde\Sigma$ may have ``necks''
which shrink off on approach to $\Sigma$, so the union of several components of $\Sigma$ may be close to a
single component of $\wtilde\Sigma$.
  
{\bf{}(2)\,} For $\Sigma\in{\cal{}E}_{0}(\Lambda)$ we let $\Omega_{\Sigma}=\cup_{j=1}^{Q}\Omega_{j}$ and we
let $\varphi_{\delta}\in W^{1,2}_{0}(\Omega_{\Sigma})$ be defined by
$$%
\varphi_{\delta}|\Omega_{j}=\varphi_{j,\delta}, 
$$%
where $\Omega_{j},\varphi_{j,\delta}$, $j=1,\ldots,Q$ are as in Lemma~4.4.  Then $\varphi_{\delta}$ satisfies
$$% 
\int_{\Omega_{\Sigma}}\Bigl(\nabla_{\Sigma} \varphi_{\delta}\cdot \nabla_{\Sigma}v - %
        \bigl(q_{\Sigma}+\lambda_{1,\delta}(\Sigma)\bigr) v\varphi_{\delta}\Bigr)\le 0,  %
                                                           \quad v\ge 0,\,\, v\in W^{1,2}_{0}(\Omega_{\Sigma}), %
\leqno{4.5}
 $$% 
 (i.e.\ weakly satisfies
$\Delta_{\Sigma}\varphi_{\delta}+(q_{\Sigma}+\lambda_{1,\delta}(\Sigma))\varphi_{\delta}\ge 0$ on
$\Omega_{\Sigma}$), where
$$%
\lambda_{1,\delta}(\Sigma)=\max\{\lambda_{1}(\Omega_{j}):j=1,\ldots,Q\} 
\leqno{4.6}
$$%
with $\lambda_{1}(\Omega_{j})$ the first eigenvalue of $\Omega_{j}$ as in Lemma~4.4. We claim that in fact
 $$%
 \int_{\Omega_{\Sigma}}\Bigl(\nabla_{\Sigma} \varphi_{\delta}\cdot \nabla_{\Sigma}v - %
                                      \bigl(q_{\Sigma}+\lambda_{1,\delta}(\Sigma)\bigr) v\varphi_{\delta}\Bigr)\le 0  %
\leqno{4.7}
 $$% 
 for any non-negative $C^{1}(\overline\Omega_{\Sigma})$ function $v$ (without the assumption that $v$
vanishes on $\partial\Omega_{\Sigma}$).  We can check this by replacing $v$ in~4.5 by
$v{}f_{k}\circ\varphi_{\delta}(\in{}W^{1,2}_{0}(\Omega_{\Sigma}))$, where $f_{k}$ is the piecewise linear
function $f_{k}(t) =\max\{0,\min\{k(t-\tfrac{1}{k}),1\}\}$.  Since
$$%
\int_{\{x\in\Omega_{\Sigma}:\varphi_{\delta}(x)<\tfrac{2}{k}\}}  %
                |\nabla\varphi_{\delta}\cdot \nabla v|  %
       \le C\|\varphi_{\delta}\|_{W^{1,2}(\Omega_{\Sigma})} %
                   ({\cal{}H}^{n}\{x\in\Omega_{\Sigma}:\varphi_{\delta}(x)<\tfrac{2}{k}\})^{1/2}\to 0  %
$$%
and $\nabla_{\Sigma}\varphi_{\delta}\cdot\nabla_{\Sigma}(vf_{k}\circ\varphi_{\delta})=
f_{k}\circ\varphi_{\delta}\nabla_{\Sigma}\varphi_{\delta}\cdot \nabla_{\Sigma}v+
vf_{k}^{\prime}(\varphi_{\delta})|\nabla_{\Sigma}\varphi_{\delta}|^{2}\ge
f_{k}\circ\varphi_{\delta}\nabla_{\Sigma}\varphi_{\delta}\cdot \nabla_{\Sigma}v$, we see that~4.5 gives 4.7 as
claimed in the limit as $k\to \infty$.

Note also that in accordance with the conclusions of Lemma~4.4 for this $\delta>0$ we have
$$%
\left\{
\begin{aligned}%
&\max_{\Omega_{\Sigma}}\varphi_{\delta}=1\text{ and  }  %
     -(\Delta_{\Sigma}\varphi_{\delta}+ q_{\Sigma}\varphi_{\delta})=  %
                  \lambda_{1,\delta}(\Sigma)\varphi_{\delta}% 
    \text{ weakly on }\Omega_{\Sigma}, \\ %
&{\cal{}H}^{n}(\Omega_{\Sigma}\setminus\{x\in\Omega_{\Sigma}:\varphi_{\delta}(x)>\tau\}) <\delta %
\end{aligned}\right.%
\leqno{4.8}
$$%
for suitable $\tau=\tau(\delta,\Lambda,{\cal{}E}_{0},{\cal{}Q}_{0})>0$ and all $\Sigma\in {\cal{}E}_{0}(\Lambda)$.

\medskip

There are also various circumstances which make it possible to prove upper bounds for $\lambda_{1}(\Sigma)$
to complement the lower bound~4.3. For example, if $\lambda_{0}\in \R$, $\C\in{\cal{}C}_{0}$, $v \in
W^{1,2}_{\text{loc}}(\Sigma)\cap L^{4}(\Sigma)\setminus\{0\}$, $q_{\Sigma}v^{2}\in L^{1}(\Sigma)$, and if $v$ is
a subsolution of the equation $\Delta_{\Sigma}u+(q_{\Sigma}+\lambda_{0})u=0$ in the sense that
$\int_{\Sigma}(\nabla v\cdot \nabla \zeta- (q_{\Sigma}+\lambda_{0})v\zeta)\le 0$, whenever $\zeta$ is a bounded
non-negative $W^{1,2}$ function with compact support in $\Sigma$, and if ${\cal{}H}^{n-4}(\sing\C)=0$, then
$\lambda_{1}(\Sigma)\le \lambda_{0}$.  To see this, we let $s_{\delta}$ be a function analogous to that used in
the above discussion, except that we now choose the balls $B_{\rho_{j}/2}(y_{j})$ to cover $\sing \C$ and
$\sum_{j}\rho_{j}^{n-4}<\delta$.  Then using the above inequality with $v\, s^{2}_{\delta}$ in place of $\zeta$,
we infer
$$% 
\int_{\Sigma} s_{\delta}^{2}\bigl(|\nabla v|^{2}- (q_{\Sigma}+\lambda_{0})\, v^{2}\bigr)%
\le -\int_{\Sigma}( 2 vs_{\delta}\nabla v\cdot \nabla s_{\delta}) %
$$% 
and hence by the Cauchy-Schwarz inequality we infer first that
\begin{align*} %
  \int_{\Sigma}\bigl({1\over{}2} s_{\delta}^{2}|\nabla v|^{2}- %
  (q_{\Sigma}+\lambda_{0})\, v^{2}\bigr)%
  & \le C\int_{\Sigma} |\nabla s_{\delta}|^{2}v^{2}\\ %
  &\le C\Bigl(\int_{\Sigma}|\nabla s_{\delta}|^{4} \Bigr)^{1/2}\Bigl(\int_{\Sigma}v^{4}\Bigr)^{1/2}. %
\end{align*} %
Since the right side $\to 0$ as $\delta\downarrow 0$ this then implies by Fatou's lemma that
$\int_{\Sigma}|\nabla v|^{2}<\infty$, and then going back to the first inequality above we have
\begin{align*} %
  \int_{\Sigma} s_{\delta}^{2}\bigl(|\nabla v|^{2}- %
                                                             (q_{\Sigma}+\lambda_{0})\, v^{2}\bigr)%
  &\le \int_{\Sigma}\bigl( 2 vs_{\delta}\nabla v\cdot \nabla s_{\delta}\bigr)\\ %
  &\le C\Bigl(\int_{\Sigma} |\nabla v|^{2}\Bigr)^{1/2}\Bigl(\int_{\Sigma} |\nabla s_{\delta}|^{2}\Bigr)^{1/2} \to 0  %
                                                                                                 \text{ as }\delta\downarrow 0.%
\end{align*} %
On the other hand $\int_{\Sigma}|\nabla (s_{\delta}v)|^{2}= \int_{\Sigma}(s_{\delta}^{2}|\nabla
v|^{2}+2s_{\delta}v\nabla s_{\delta}\cdot \nabla v+ v^{2}|\nabla s_{\delta}|^{2})$ and hence the above inequality
implies
$$% 
 \limsup_{\delta\downarrow 0}\int_{\Sigma}\bigl(|\nabla (s_{\delta}v)|^{2}  %
                                                        -(q_{\Sigma}+\lambda_{0})(s_{\delta}v)^{2}\bigr) \le 0, %
$$%
 whence, since $\int_{\Sigma}(s_{\delta}v)^{2}\to \int_{\Sigma}v^{2}\in (0,\infty)$,
$$%
 \limsup_{\delta\downarrow0}\Bigl(\int_{\Sigma}(s_{\delta}v)^{2}\Bigr)^{-1}   %
         \int_{\Sigma} \bigl(|\nabla(s_{\delta} v)|^{2}- q_{\Sigma} (s_{\delta}v)^{2}\bigr)\le \lambda_{0}. %
$$%
 Thus we have proved that
$$% 
\lambda_{1}(\Sigma) \le \lambda_{0}  %
\leqno{4.9} %
$$%
if ${\cal{}H}^{n-4}(\sing \C)=0$ and if there is a non-negative $W^{1,2}_{\text{loc}}(\Sigma)\cap
L^{4}(\Sigma)\setminus\{0\}$ subsolution of the equation $\Delta_{\Sigma}u+r^{-2}(q_{\Sigma}+\lambda_{0})u=0$
on $\C$.

\section{A partial Harnack theory}

Here $M$ will denote any fixed element of ${\cal{}P}$ and $a,b_{1},\ldots,b_{n}$ will be measurable functions
on~$M$ with $|a|^{1/2}+|b|$ locally integrable on $M$, $b=(b_{1},\ldots,b_{n})$. We also need to assume
$L^{n+\alpha}$ ($\alpha>0$) bounds on the function $|a|^{1/2}+|b|$ on the various domains which arise
here---the precise bounds needed will be stipulated in each case.

Recall that $u$ is a positive supersolution of the equation $\Delta_{\! M}u + b\cdot\nabla_{\!M}u+au=0$ on
$M$ means that $u>0$ a.e.\ on $M$, $u\in W^{1,2}_{\text{loc}}(M)$ and
$$% 
\int_{\! M}\bigl(a u\, \zeta +b\cdot \nabla_{\! M}u\zeta- \nabla u \cdot \nabla \zeta\bigr)\le 0 %
\leqno{5.1} %
$$%   
for each non-negative $C_{c}^{1}(M)$ function $\zeta$.

Of course by classical De~Giorgi Nash theory (applied locally on the $C^{1}$ manifold $M$) there is no loss of
generality in assuming that $u$ is positive on $M$ with local uniform positive lower bounds.  Also if $p\in (0,
1)$ and if we use $\zeta u^{p-1}$ in place of $\zeta$ in this inequality, then, again letting $\epsilon\downarrow
0$, we get the inequality
$$% 
\int_{\! M}\zeta\bigl(p au^{p}+ b\cdot \nabla_{\! M}u^{p} + {4(1-p)\over p}|\nabla u^{p/2}|^{2}\bigr)\le  %
                               \int_{\! M}\nabla u^{p}\cdot \nabla \zeta, %
\leqno{5.2} %
$$% 
for each non-negative $\zeta\in C_{c}^{1}(M)$.  Now $\nabla u^{p}= 2 u^{p/2}\nabla u^{p/2}$, so if we replace
$\zeta$ by $\zeta^{2}$ and use the Cauchy-Schwarz inequality then we obtain
$$%    
\int_{\! M}\zeta^{2}\bigl( {(1-p)\over p}|\nabla u^{p/2}|^{2}\bigr)\le   %
                      8(1-p)^{-1}\int_{\! M}u^{p}\bigl(|\nabla \zeta|^{2}+(|a|+|b|^{2})\zeta^{2}\bigr), %
\leqno{5.3} %
$$% 
for each $\zeta\in C_{c}^{1}(M)$.  There is also a version of this for $p=0$:
$$%    
\int_{\! M}\zeta^{2}\bigl( |\nabla w|^{2}\bigr)\le  %
    8\int_{\! M}\bigl(|\nabla \zeta|^{2}+(|a|+|b|^{2})\zeta^{2}\bigr),  %
\leqno{5.4} %
$$%
where $w=\log u$, which is obtained by substituting $\zeta\, u^{-1}$ in place of $\zeta$ in~5.1.

Now we claim that~5.4 is valid for $\zeta\in C^{1}_{c}(U_{\! M})$.  That is, it is not necessary
that $\zeta$ vanish in a neighborhood of $\sing M$. Indeed to see that~5.4 is valid we simply replace $\zeta$
by $\zeta\, s_{\delta}$ and let $\delta\downarrow 0$ as in our discussion of~3.2 

We cannot do quite the same thing to justify the fact that~5.3 holds for any $\zeta\in C^{1}_{c}(U_{\! M})$,
because $u^{p}$ is not necessarily bounded. However notice that if we take any $K>0$ and any $C^{2}$ concave
increasing function $f_{K}$ with $f_{K}(t)\equiv t$ for $0\le t\le K$ and $f_{K}(t)\equiv K+{1\over 2}$ for
$t>K+1$, then by replacing $\zeta$ by $\zeta\,f_{K}^{\prime}(u)$ in~5.1 we obtain directly that $f_{K}(u)$ is a
supersolution of an equation of the same form.  Thus we have the estimate~5.4 with $u_{K}=f_{K}(u)$ in place
of $u$ and with any $\zeta\in C^{1}_{c}(M)$:
$$% 
\int_{\! M}\zeta^{2}( |\nabla u_{K}^{p/2}|^{2})\le C\int_{\! M}u_{K}^{p}\bigl(|\nabla
\zeta|^{2}+(|a|+|b|^{2})\zeta^{2}\bigr),
$$% 
for any $\zeta\in C^{1}_{c}(M)$.  Now with $s_{\delta}$ as in the discussion following~3.2, for any $z\in
C^{1}_{c}(U_{\! M})$ we can then substitute $\zeta\, s_{\delta}$ in place of $\zeta$ here.  Since $u_{K}$ is
bounded, we can then let $\delta\downarrow 0$ (as we did to justify~3.2) to deduce that~5.3 holds for all
$\zeta\in C_{c}^{1}(U_{\! M})$, with $u_{K}$ in place of $u$. Then letting $K\uparrow \infty$, we deduce
that~5.3 also holds for all $\zeta\in C_{c}^{1}(U_{\!  M})$ provided $u^{p}\in L^{1}(M\cap K)$ for each compact
$K\subset U_{M}$. We show below in~5.7 that indeed $u^{p}\in L^{1}(M\cap K)$ for each $p\in
(0,{n\over{}n-2})$.

\smallskip

{\bf 5.5 Remark:} Notice that in particular~5.2 says that $u^{p}$ is a supersolution of the equation $\Delta_{\! M}
u +b\cdot \nabla_{\! M}u+ p a u= 0$ locally in $M$ for each $p\in (0,1)$.

\smallskip

Now assume that $y\in\overline M\cap U_{M}$,  take any closed ball $B_{\rho}(y)\subset U_{M}$ and assume
$$%
\rho^{-n/(n+\alpha)}\| |a|^{1/2}+|b|\|_{L^{n+\alpha}(M\cap B_{\rho}(y))}\le \beta %
\leqno{5.6}
$$%
for some constants $\alpha,\beta>0$.  Recall that the weak Harnack theory for supersolutions on domains in
$\R^{n}$ says that, for any $\theta\in ( 0,1)$, $\inf_{B_{\theta\rho}(y)}v\ge
C^{-1}(\rho^{-n}\int_{B_{\theta\rho}}v^{p})^{1/p}$ for positive supersolutions $v$ of the equation $\Delta
u+b\cdot Du + au=0$, with $C=C(n,\theta,p,\beta)$, assuming $a,b$ satisfy an inequality like~5.6 in the
Euclidean ball $B_{\rho}(y)$ (rather than in $M\cap B_{\rho}(y)$). The proof of this requires not only a Sobolev
inequality (analogous to the inequality established for surfaces $P\in{\cal{}P}$ in~2.6$^{\prime}$) but also a
Poincar\'e inequality, and unfortunately in the present setting there is no such inequality (this would require
strong connectivity hypotheses on the submanifolds in the class ${\cal{}P}$), so we cannot follow the $\R^{n}$
procedure to give a Harnack theory.  Nevertheless, with only the Sobolev inequality of Remark~2.6$^{\prime}$
and the modified Poincar\'e inequality of Remark~2.6 at our disposal, we claim that it is still possible to prove
the following:

\begin{state}{\bf{}5.7 Lemma.} %
  Suppose ${\cal{}P}$ is a regular multiplicity~1 class (so that 1.1--1.4 hold), $\theta\in (0,1)$, $p\in
[1,{n\over{}n-2})$, $M\in{\cal{}P}$, $y\in \overline M$ with $B_{\rho}(y)\subset U_{M}$, and $u$ is a positive
supersolution of the equation $\Delta_{M}u+b\cdot \nabla_{M}u+ au=0$ on $M\cap \breve B_{\rho}(y)$, where
$a,b$ satisfy~5.6. Then there is $\delta=\delta({\cal{}P},\alpha,\beta,\theta,p)\in (0,\ha]$ such that 
$$% 
\rho^{-n}\int_{\! M\cap B_{\theta\rho}(y)} u^{p} \le C\lambda^{p}, \quad C=C({\cal{}P},\alpha,\beta,\theta,p), %
$$% %
for any $\lambda$ such that ${\cal{}H}^{n}(\{x\in M\cap B_{\rho}(y): u(x)>\lambda\})\le\delta\rho^{n}$.
\end{state}%

{\bf{}Proof:} We use some modications of the relevant part of the De~Giorgi Nash Moser theory
(\cite[\S8.8--\S8.10]{GilT83}).  First, by rescaling, we can assume $\rho=1$, so we aim to prove a bound for the
$L^{p}$ norm of $u$ over the ball $B_{\theta}(y)$, where $\theta\in (0,1)$.  For $\lambda>0$ (fixed for the
moment), let
$$%
\wtilde w= \min\{\max\{\log(u/\lambda), 0\},K\}, %
$$%
where $K\ge 2$ (we plan to let $K\uparrow\infty$ eventually) and observe that for any $q\ge 1$, $ m,
\gamma>0$ with $2 m -\gamma> 2$, and non-negative $\zeta\in C^{1}_{\text{c}}(\breve B_{\theta}(y))$ satisfying
$$%
\kappa\gamma=\gamma+2,\quad \text{$\zeta\equiv 1$ on $B_{\theta^{2}}$, and $0\le \zeta\le 1$, $|\nabla
\zeta|\le C(\theta)$ on $B_{1}$,} %
\leqno{(1)}
$$%
where $\kappa=n/(n-2)$ for $n\ge 3$ as in~2.5$^{\prime}$, 2.6$^{\prime}$ and $\kappa>2$ is arbitrary in case
$n=2$. (So $\gamma=n-2$ if $n\ge 3$ and $\gamma=2/(\kappa-1)$ if $n=2$.)  For the remainder of the proof we
let $C$ denote any constant 
$$%
C=C({\cal{}P},\alpha,\beta,m,N,\theta,p);
$$%
it is important to keep track of the $q$ dependence though, so that will be indicated explicitly at each stage of
the proof. Then we can
apply~2.6$^{\prime}$ with $\varphi=\wtilde w^{2q}\zeta^{2 m q-\gamma}$ with $q\ge 1$, giving
$$%
\bigl(\tint\, (\wtilde w\zeta^{ m})^{2\kappa q}\,d\mu\bigr)^{1/\kappa} \le  %
                    C q^{2} \tint\, \bigl(|\nabla \wtilde w|^{2} \wtilde w^{2q-2}\zeta^{2q m+2} +  %
                                     \wtilde w^{2q}\zeta^{2q m}\bigr)\,d\mu,    %
\leqno{(2)}
$$%
where $\mu$ is the Borel measure defined by
$$%
\mu=\zeta^{-\kappa\gamma}{\cal{}H}^{n}\res (M\cap B_{1}(y))  %
                             (=\zeta^{-\gamma-2}{\cal{}H}^{n}\res (M\cap B_{1}(y))). %
$$%
To handle the first term on the right of~(2) we proceed slightly differently in the cases $q\ge 2$  and $q\in
[1,2)$.  If $q\ge 2$, by replacing $\zeta$ in~5.4 with $\wtilde w^{q-1}\zeta^{q m-\gamma/2}$
then we get
$$%
\tint (\wtilde w^{2q-2}\zeta^{2 m q+2}|\nabla \wtilde w|^{2})\,d\mu \le  %
 Cq^{2}\tint (\wtilde w^{2q-4}\zeta^{2 m q+2}|\nabla \wtilde w|^{2}+A\wtilde w^{2q-2}\zeta^{2 m
q})\,d\mu,
$$% 
where
$$%
A= 1+(|a| +|b|^{2})\zeta^{2},
$$%
and since  $\wtilde w^{2q-2}\le 1+\wtilde
w^{2q}$ this gives
$$%
\tint (\wtilde w^{2q-2}\zeta^{2 m q+2}|\nabla \wtilde w|^{2})\,d\mu \le  %
Cq^{2}+ Cq^{2}\tint (\wtilde w^{2q-4}\zeta^{2 m q+2}|\nabla \wtilde w|^{2}+A\wtilde w^{2q}\zeta^{2 m
q})\,d\mu.
\leqno{(3)}
$$% 
For $q\ge 2$ Young's inequality gives
$$%
Cq^{2}\wtilde w^{2q-4} \le \ha \wtilde w^{2q-2} + \wtilde C^{q}q^{2q},
$$%
and we thus get for any $q\ge 2$
$$%
\tint (\wtilde w^{2q-2}\zeta^{2 m q+2}|\nabla \wtilde w|^{2})\,d\mu   %
 \le Cq^{2}+C^{q}q^{2q}\tint (\zeta^{2 m q+2}|\nabla \wtilde w|^{2})\,d\mu+  %
                                         Cq^{2}\tint (A \wtilde w^{2q}\zeta^{2 m q})\,d\mu.  %
\leqno{(4)}
$$%
On the other hand if $q\in [1,2)$ then we can first use $\wtilde w^{2q-2}\le 1+\wtilde w^{2}$ and hence
$$%
\tint (\wtilde w^{2q-2}\zeta^{2 m q+2}|\nabla \wtilde w|^{2})\,d\mu \le  
\tint (\wtilde w^{2}\zeta^{2 m q+2}|\nabla \wtilde w|^{2})\,d\mu +\tint (\zeta^{2 m q+2}|\nabla \wtilde
w|^{2})\,d\mu,  %
$$%
and by replacing $\zeta$ in~5.4 with $\wtilde w\zeta^{qm-\gamma/2}$ we obtain 
$$%
\tint (\wtilde w^{2}\zeta^{2 m q+2}|\nabla \wtilde w|^{2})\,d\mu \le  %
 C\tint (\zeta^{2 m q+2}|\nabla \wtilde w|^{2}+A\wtilde w^{2}\zeta^{2 m
q})\,d\mu,
$$% 
and hence since $\wtilde w^{2}\le 1+w^{2q}$ we again get~(4), so in fact~(4) is valid for any $q\ge 1$.

Another application of~5.4 (this time with $\zeta^{ m q-\gamma/2}$ in place of $\zeta$) gives
$$%
\tint (\zeta^{2 m q+2}|\nabla \wtilde w|^{2})\,d\mu \le Cq^{2} \tint A \zeta^{2 m q}\,d\mu \le Cq^{2},
$$%
so~(4) implies
$$%
\tint (\wtilde w^{2q-2}\zeta^{2 m q+2}|\nabla \wtilde w|^{2})\,d\mu \le  %
           C^{q}q^{2q}    +Cq^{2}\tint A (\wtilde w\zeta^{ m})^{2q}\,d\mu %
\leqno{(5)}
$$% 
for all $q\ge 1$.

Combining~(5)  with~(2) we obtain
$$%
\bigl(\tint\, (\wtilde w\zeta^{ m})^{2\kappa q}\,d\mu\bigr)^{1/\kappa} \le %
                              C^{q}q^{2q} +Cq^{2}\tint A (\wtilde w\zeta^{ m})^{2q}\,d\mu.  %
$$%
Using the H\"older inequality and 5.6 again we then conclude that 
$$%
\Bigl(\tint\, (\wtilde w\zeta^{ m})^{2\kappa q}\,d\mu\Bigr)^{1/\kappa} \le %
       C^{q}q^{2q} +Cq^{2}\tint (\wtilde w\zeta^{ m})^{2q}\,d\mu+  %
           Cq^{2}\Bigl(\tint (\wtilde w\zeta^{ m})^{2q\lambda}\,d\mu\Bigr)^{1/\lambda} %
$$%
for some $\lambda=\lambda(\alpha,n)\in (1,\kappa)$.  Using the
$L^{p}$ interpolation inequality (with respect to the measure $\mu$)
$$%
\|f\|_{L^{\lambda}(\mu)} \le \epsilon \|f\|_{L^{\kappa}(\mu)} + \epsilon^{-C(\kappa,\lambda)} \|f\|_{L^{1}(\mu)} 
$$%
with $f=(\wtilde w \zeta^{ m})^{2q}$ then yields
$$%
\Bigl(\tint (\wtilde w\zeta^{ m})^{2q\kappa}\,d\mu\Bigr)^{1/\kappa} \le  C^{q}q^{2q} + C
q^{C}\tint(\wtilde w\zeta^{ m})^{2q}\,d\mu,
$$%
and hence
$$%
\Psi(\kappa q) \le C q  + C^{1/q}q^{C/q} \Psi(q)\text{ with } %
                                        \Psi(q)=\Bigl(\tint(\wtilde w\zeta^{ m})^{2q}\,d\mu\Bigr)^{1/2q},  %
\leqno{(6)}
$$%
for each $q\ge 1$.  Replacing $q$ by $\kappa^{\nu},\,\nu=0,1,2,\ldots$, we obtain
$$%
\Phi(\nu+1)\le C\kappa^{\nu}+C^{1/\kappa^{\nu}}\kappa^{\nu/\kappa^{\nu}}\Phi(\nu),\,\, %
                              \nu=0,1,2,\ldots, \text{ with } \Phi(\nu)= \Psi(\kappa^{\nu}). %
$$%
Iterating, and using the facts that $\sum_{j=0}^{\nu}\kappa^{j}\le C\kappa^{\nu}$ and
$\sum_{}\nu\kappa^{-\nu}<\infty$, we get
$$%
\Phi(\nu)  \le C\kappa^{\nu}+C\Phi(0).
$$%
But 
$$%
\Phi(0)=\tint(\wtilde w\zeta^{ m})^{2}\,d\mu\le\int_{M\cap B_{\theta}(y)}\wtilde w^{2}\,d{\cal{}H}^{n}\le  %
                       C\int_{M\cap B_{(1+\theta)/2}}|\nabla \wtilde w|^{2}\,d{\cal{}H}^{n}\le C  %
$$%
by~2.6(2), 5.4, and~5.6, so in fact
$$%
\Phi(\nu)  \le C\kappa^{\nu}, \quad \nu=0,1,\ldots,
$$%
and then in view of~(6)  we conclude
$$%
\tint(\wtilde w\zeta^{ m})^{2q}\,d\mu  \le C^{q}q^{2q}, \quad q=0,1,\ldots,
$$%
and since $\zeta\equiv 1$ on $B_{\theta^{2}}(y)$, this gives
$$%
\int_{M\cap B_{\theta^{2}}(y)} \wtilde w^{2q}\,d{\cal{}H}^{n}  \le C^{q}q^{2q}, \quad q=0,1,\ldots.
$$%
Since $q^{2q}\le C^{q}(2q)!$ we can sum over $q$ here to conclude that
$$%
\int_{M\cap B_{\theta^{2}}(y)} e^{p_{0}\wtilde w}\,d{\cal{}H}^{n}  \le C
$$%
for some $p_{0}=p_{0}({\cal{}P},\alpha,\beta,\theta)\in (0,1)$, and since $\wtilde w\uparrow 
\max\{\log({u\over{}\lambda}),0\}$ as $K\to \infty$, this implies
$$%
\int_{M\cap B_{\theta^{2}}(y)} u^{p_{0}}\,d{\cal{}H}^{n}  \le C \lambda^{p_{0}}.
\leqno{(7)}
$$%
However~5.3 and the Sobolev inequality of Remark~2.6$^{\prime}$ evidently imply that
$$%
\Bigl(\int_{M\cap B_{\theta^{2}}(y)}u^{\kappa p}\Bigr)^{1/\kappa p} \le  %
                                    C\Bigl(\int_{M\cap B_{\theta}(y)}u^{p}\Bigr)^{1/p},  \quad p\in (0,1),%
$$%
where $C=C({\cal{}P},\alpha,\theta,\beta,p)$, and by finite iteration of this inequality we have 
$$%
\Bigl(\int_{M\cap B_{\theta^{j+1}}(y)}u^{\kappa^{j} p}\Bigr)^{1/\kappa^{j} p} \le  %
      C\Bigl(\int_{M\cap B_{\theta}(y)}u^{p}\Bigr)^{1/p},  \quad \text{ provided }\kappa^{j-1} p\in (0,1),%
\leqno{(8)}
$$%
where $C({\cal{}P},\alpha, j,p,\theta,\beta)$, which when used in combination with~(7) evidently implies
$$%
\int_{M\cap B_{\theta}(y)} u^{p}\,d{\cal{}H}^{n}  \le C \lambda^{p}
$$%
for each $p\in (0,\kappa)$ and each $\theta\in (0,1)$, where $C=C({\cal{}P},p,\alpha,\beta,\theta)$, thus
completing the proof of~5.7.

\medskip

We now want to show that~5.7 eliminates the possibility of concentration of $L^{p}$ norm in regions of small 
measure. Specifically, we have the following corollary. 

\begin{state}{\bf 5.8 Corollary.} % 
  If the hypotheses are as in~5.7, then there are constants $C=C(\theta, p, {\cal{}P},\alpha,\beta)>0$ and
$\delta=\delta({\cal{}P},\theta,p,\alpha,\beta)\in (0, 1)$ such that 
$$% 
\|u\|_{L^{p}(M\cap B_{\theta\rho}(y))} \le C \|u\|_{L^{1}(M\cap B_{\rho}(y)\setminus \Omega_{\delta})}
$$% 
whenever $\Omega_{\delta}\subset M\cap B_{\rho}(y)$ has ${\cal{}H}^{n}$-measure less than $\delta\rho^{n}$.
\end{state}

{\bf{}Proof:} By change of variable $x\mapsto \eta_{y,\rho}(x)$ we reduce to the case $y=0,\rho=1$, so $0\in
\overline M\cap U_{\! M}$ and $U_{\! M}\supset B_{1}$. If $\Omega_{\delta}$ is any subset of $M\cap B_{1}$ of
${\cal{}H}^{n}$-measure less than $\delta/2$, where $\delta$ is as in~2.5, and if for $K>1$ we let $A_{K}=\{x\in
M\cap B_{1}\setminus \Omega_{\delta}: u(x)> K\|u\|_{L^{1}(M\cap B_{1}\setminus \Omega_{\delta})}\}$, then
$$%
\|u\|_{L^{1}(M\cap B_{1}\setminus\Omega_{\delta})} \ge \|u\|_{L^{1}(A_{K})} \ge
K {\cal{}H}^{n}(A_{K}) \|u\|_{L^{1}(M\cap B_{1}\setminus \Omega_{\delta})},
$$%
so that ${\cal{}H}^{n}(A_{K})\le K^{-1}$ and hence with the choice $K=2/\delta$ we get
${\cal{}H}^{n}(A_{K})<\delta/2$. Thus with this $K$ we have
$$%
\begin{aligned}%
&{\cal{}H}^{n}(\{x\in M\cap B_{1}:u(x)> K\|u\|_{L^{1}(M\cap B_{1}\setminus \Omega_{\delta})}\})  \\  %
&\hskip1in \le {\cal{}H}^{n}(\{x\in M\cap B_{1}\setminus\Omega_{\delta}:u(x)> K\|u\|_{L^{1}(M\cap B_{1} %
                              \setminus \Omega_{\delta})}\}) + {\cal{}H}^{n}(\Omega_{\delta}) \\ % 
&\hskip1in ={\cal{}H}^{n}(A_{K})+{\cal{}H}^{n}(\Omega_{\delta}) \le \delta/2+\delta/2=\delta, %
\end{aligned}% 
$$%
whence we can apply~5.7 with $\lambda=K\|u\|_{L^{1}(M\cap{}B_{1}\setminus\Omega_{\delta})}$, and this gives
the required result.

\medskip

\section{Proof of Theorem 1}

Suppose the theorem is false for some given $\alpha,\beta>0$, $p\in [1,{n\over{}n-2})$, classes ${\cal{}P}$,
${\cal{}C}_{0}\subset{\cal{}C}$ and $\gamma<\gamma_{0}$. Then for each choice of $\rho\in (0,{1\over{}4})$ the
theorem fails, so there is $\tau_{k}\downarrow 0$, $M_{k}\in{\cal{}P}$ with $U_{M_{k}}\supset \breve
B_{3/2}\setminus B_{\tau_{k}}$, $\C_{k}\in {\cal{}C}_{0}$ such that 1.1--1.11 all hold with
$M_{k},\C_{k},u_{k},q_{k},a_{k},b_{k},\tau_{k}$ in place of $M,\C,u,q,a,b,\tau$ respectively, yet such that
$$%
\|u_{k}\|_{L^{p}(M_{k}\cap B_{\rho}\setminus B_{\rho/2})} <   %
                     \rho^{-\gamma} \|u_{k}\|_{L^{p}(M_{k}\cap B_{1}\setminus B_{1/2})}.  %
\leqno{(1)}
$$%
Thus
$$%
d(M_{k}\cap B_{3/2}\setminus B_{\tau_{k}},\C_{k}\cap B_{3/2}\setminus B_{\tau_{k}})\to 0
$$%
and by compactness of ${\cal{}E}_{0}$ we can pass to a subsequence and select $\C\in{\cal{}C}_{0}$
with
$$%
d(\C_{k}\cap B_{3/2},\C\cap B_{3/2})\to 0.
$$%
Let
$$%
\Sigma = \C\cap S^{N-1}.
$$%
Since $\Sigma$ has only finitely many connected components and the convergence of $M_{k}$ to $\C$ is in the
$C^{1}$ sense of~1.4 near compact subsets of $\C\cap \breve B_{3/2}$ we can use the classical Harnack theory
of De~Giorgi Nash Moser (i.e.\ \cite[\S8.8--\S8.10]{GilT83}) applied locally to the solutions $u_{k}$ on $M_{k}$,
and the local $W^{1,2}$ estimates of~5.3 for $u^{p}$ ($p<1$), together with the Rellich compactness theorem,
to assert that, for sufficiently small $\delta_{0}=\delta_{0}(\Sigma)>0$, a subsequence of the normalized sequence
$$%
\wtilde u_{k} =    %
    \|u_{k}\|^{-1}_{L^{1}(\{x\in M_{k}\cap B_{1}\setminus B_{1/2}:\dist(x,\sing \C)\ge \delta_{0}\})} u_{k}  %
$$%
converges (in the sense discussed in the remark following~1.4) strongly in $L^{p}$ for $p<{n\over{}n-2}$, on
compact subsets of $\C\cap \breve B_{3/2}$ to a non-negative $u\in W^{1,2}_{\text{loc}}(\C\cap \breve B_{3/2})$
which satisfies $u>0$ on at least one connected component of $\C$, and
$$%
\Delta_{\C} u +  r^{-2}q_{\Sigma}u \le 0  %
\leqno{(2)}
$$% 
weakly on $\C\cap \breve B_{3/2}$, where $r^{-2}q_{\Sigma}$ is the uniform limit of $q_{k}$ on compact
subsets of $\C$. Furthermore, by the inequality on the right of~2.1 and the fact that ${\cal{}H}^{n-2}(\sing \C\cap
\breve B_{3/2})<\infty$ there is $\delta(\tau)\downarrow 0$ as $\tau\downarrow 0$, with $\delta(\tau)$ not
depending on $k$, such that
$$%
{\cal{}H}^{n}(M_{k}\cap \{x\in B_{3/2}\setminus B_{\rho/4}:\dist(x,\sing\C)\ge \tau\})<\delta(\tau). 
$$%
Then the partial Harnack theory (in particular Corollary~5.8) is applicable, ensuring that in fact we have the
$L^{1}$ norm convergence
$$%
\|\wtilde u_{k}\|_{L^{1}(M_{k}\cap B_{R}\setminus B_{R/2})}\to \|u\|_{L^{1}(\C\cap B_{R}\setminus B_{R/2})}
\leqno{(3)}
$$%  
for each $R\in ({\rho\over{}2},{3\over{}2})$ (and in particular this holds with $R=1$ and $R=\rho$).

Observe that in view of~(2) we can apply~4.3 with each $\tau>0$ and so 
$$%
\lambda_{1}({\cal{}E}_{0})\ge \lambda_{1}(\Sigma)\ge -\bigl({n-2\over{}2}\bigr)^{2} %
$$%
and hence $\gamma_{0}$ (in the statement of Theorem~1) is the smaller root of the quadratic equation
$t^{2}-(n-2)t-\lambda_{1}({\cal{}E}_{0})$. Thus if we take
$$%
\mu=\ha(\gamma+\gamma_{0}), %
\leqno{(4)}
$$%
then we have 
$$%
\gamma+C<\mu<\gamma_{0}-C \text{ and }\mu^{2}-(n-2)\mu-\lambda_{1}({\cal{}E}_{0}) \ge C,
\leqno{(5)} 
$$%
with $C=C(n,\gamma,{\cal{}E}_{0})>0$.  

Now let $\varphi_{\delta},\lambda_{1,\delta}(\Sigma)$, and $\Omega_{\Sigma}\subset\Sigma$ be as in~4.7. 
Notice that the weak form of~(2) on $\C\cap\breve B_{3/2}$ is
$$%   
\int_{0}^{3/2}\int_{\Sigma} \Bigl(-u^{\prime}\zeta^{\prime}-  %
               r^{-2}\nabla_{\Sigma} u\cdot \nabla_{\Sigma} \zeta +  %
                              r^{-2}q_{\Sigma} u \zeta\Bigr)r^{n-1}\,drd\omega\le 0 %
\leqno{(6)}
$$% 
for any $\zeta\in C^{\infty}_{c}(\C\cap\breve B_{3/2})$ with $\zeta\ge 0$, where $u^{\prime}$ means ${\partial
u\over \partial r}$.  Replacing $\zeta$ by $\zeta_{1}(r)\varphi_{\delta}$ in~(6) gives
$$%
\int_{0}^{3/2}\int_{\Sigma} \Bigl(-u^{\prime}\zeta_{1}^{\prime}\varphi_{\delta} - 
                             r^{-2}\zeta_{1}\nabla_{\Sigma}u\cdot \nabla_{\Sigma} \varphi_{\delta}%
                                + r^{-2}q_{\Sigma}\zeta_{1}u \varphi_{\delta}\Bigr)r^{n-1}\,drd\omega   \le 0. %
$$%
Using inequality~4.7 with $v(\omega)=u(r\omega)$ on the left, and writing $v_{\delta}(r)=\langle u(r),
\varphi_{\delta}\rangle_{L^{2}(\Omega_{\Sigma})}$, we then conclude
$$% 
-\int_{0}^{3/2}\zeta_{1}^{\prime}v_{\delta}^{\prime}r^{n-1}\,dr -  %
               \lambda_{1,\delta}(\Sigma)\int_{0}^{3/2}\zeta_{1} v_{\delta}r^{n-3}\,dr \le 0. %
$$% 
That is, weakly $v_{\delta}$ satisfies
$$% 
r^{1-n}(r^{n-1}v_{\delta}^{\prime})^{\prime} -  %
          \lambda_{1,\delta}(\Sigma) r^{-2} v_{\delta} \le 0,\quad r\in (0,3/2)% 
\leqno{(7)} %
$$% 
Now with $\mu$ as in~(5) let
$w_{\delta}=r^{\mu}v_{\delta}$. Then  $v_{\delta}=r^{-\mu}w_{\delta}$ and so
$$%
 v_{\delta}^{\prime} = -\mu r^{-\mu-1}w_{\delta}+ r^{\mu} w_{\delta}^{\prime},\quad  %
      v_{\delta}^{\prime\prime}= \mu(\mu+1)r^{-\mu-2}w_{\delta}-2\mu  %
                               r^{-\mu-1}w_{\delta}^{\prime}+ r^{-\mu}w_{\delta}^{\prime\prime}. %
$$%
 Then substituting in the differential inequality~(7) we get
$$% 
 r^{-2}(\mu^{2}-(n-2)\mu - \lambda_{1,\delta}(\Sigma))w_{\delta} +   %
            (n-1-2\mu)r^{-1}w_{\delta}^{\prime} + w_{\delta}^{\prime\prime} \le 0.  %
$$% 
Since $\mu^{2}-(n-2)\mu - \lambda_{1,\delta}(\Sigma) > 0$ (by~(5)), we see that the previous implies
 $$% 
 r^{-(1+\beta)}(r^{1+\beta}w_{\delta}^{\prime})^{\prime}\le 0 \quad \mbox{ (weakly for $r\in (0,3/2)$) },
 $$% 
 where $\beta= n-2-2\mu> 0$.  This says that $w_{\delta}$ is a concave function of the new variable
$s=r^{-\beta}\in ((3/2)^{-\beta},\infty)$, and since $w_{\delta}$ is non-negative we see that then $w_{\delta}$
must be increasing with respect to the variable $s$; that is $w_{\delta}$ is a decreasing function of the variable
$r$, so that
$$% 
w_{\delta}(r_{1})\ge w_{\delta}(r_{2})\mbox{ for }0<r_{1}<r_{2}<3/2,
$$% 
which in terms of $u$ says
$$% 
r_{1}^{\mu}\int_{\Sigma}u(r_{1}\omega)\,\varphi_{\delta}(\omega)\ge 
          r_{2}^{\mu}\int_{\Sigma}u(r_{2}\omega)\,\varphi_{\delta}(\omega)\mbox{ for }0<r_{1}<r_{2}<3/2. 
$$% 
Integrating over $(1/2, 1)$ with respect to the variable $r_{1}$ and over $(\rho/2,\rho)$ with respect to the
variable $r_{2}$, we then conclude that
$$% 
\rho^{\mu-n}\int_{B_{\rho}\setminus B_{\rho/2}}u\,\varphi_{\delta}\ge 
                                2^{-\mu-n}\int_{B_{1}\setminus B_{1/2}}u\,\varphi_{\delta}, %
$$% 
so that by~4.8 and~5.8 (applied with $M=\C$, $p=1$ and $a=b=0$) we have
$$% 
\rho^{-n}\int_{B_{\rho}\setminus B_{\rho/2}}u\ge C\rho^{-\mu}
\int_{B_{1}\setminus B_{1/2}}u,
$$% 
where $C=C({\gamma,\cal{}P}, {\cal{}E}_{0},{\cal{}Q_{0}},\alpha,\beta)$ and provided
$\delta=\delta({\gamma,\cal{}P}, {\cal{}E}_{0},{\cal{}Q_{0}},\alpha,\beta)>0$ is chosen suitably, and hence by the
norm convergence~(3) we have
$$% 
\rho^{-n}\int_{B_{\rho}\setminus B_{\rho/2}}u_{k}\ge C\rho^{-\mu} \int_{B_{1}\setminus B_{1/2}}u_{k}  %
\leqno{(8)} %
$$% 
for all sufficiently large $k$, where $C=C(\gamma,{\cal{}P},{\cal{}C}_{0},{\cal{}Q}_{0},\alpha,\beta)$,
contradicting~(1) in the case $p=1$ for sufficiently small $\rho$ (depending only on
$\gamma,{\cal{}P},{\cal{}C}_{0},{\cal{}Q}_{0},\alpha,\beta$), by~(5).  This completes the proof in the case~$p=1$.

To handle the remaining $p\in (1,{n\over{}n-2})$, observe that we could have integrated with respect to $r_{2}$
over $(1/4,5/4)$ and also we can apply the H\"older inequality on the left side of~(8), whence
$$% 
\|u_{k}\|_{L^{p}(M_{k}\cap B_{\rho}\setminus B_{\rho/2})}\ge  %
                    C\rho^{-\mu} \int_{M_{k}\cap B_{5/4}\setminus B_{1/4}}u_{k}  %
 $$% 
 for any $p\in [1,{n\over{}n-2})$.  Then by applying inequality~(8) in the proof of Lemma~5.7 (with a suitably
scaled version of $M_{k}$ in place of $M$) we conclude
$$% 
\|u_{k}\|_{L^{p}(M_{k}\cap B_{\rho}\setminus B_{\rho/2})} \ge %
                    C\rho^{-\mu} \|u_{k}\|_{L^{p}(M_{k}\cap B_{1}\setminus B_{1/2})}  % 
$$% 
for any $p=[1,{n\over{}n-2})$, with $C=C(p,\gamma,{\cal{}P},{\cal{}C}_{0},{\cal{}Q}_{0},\alpha,\beta)$,
which again contradicts~(1) for sufficiently large $k$.

\section{Application to growth estimates for exterior solutions}

In this section we want to show how the main theorem applies to give lower growth estimates for entire and
exterior solutions of the minimal surface equation. Thus we assume that $u$ is $C^{2}$ and satisfies the
minimal surface equation
$$% 
\sum_{i=1}^{n} D_{i}\Bigl( {D_{i}u\over \sqrt{1+|Du|^{2}}}\Bigr)=0 \mbox{ in }%
\R^{n}\setminus B_{1}.%
\leqno{7.1} %
$$%
We need the non-trivial general facts given in the following two theorems:

\begin{state}{\bf 7.2 Theorem.}  %
  There is a regular multiplicity~1 class ${\cal{}P}$ (as in~\S$1$) with $N=n+1$ such that ${\cal{}P}$ contains
each minimal graph $G=\graph u$, corresponding to a solution $u\in C^{2}(\Omega)$ of~7.1, where $\Omega$
is any open set in $\R^{n}$; in this case we always have that $U_{G}$ (the open set associated with
$G\in{\cal{}P}$ as in~\S$1$) is just $\Omega\times \R$.  Further the class ${\cal{}P}$ can be chosen so that the
convergence $P_{k^{\prime}}\to P$ of~1.4 is actually $C^{\infty}$ (i.e.\ $C^{k}$ for each $k$) on compact
subsets of $P$ rather than merely $C^{1}$.
\end{state}

{\bf{}Proof:} This follows from the De Giorgi theory of oriented boundaries of least area (also known as area
minimizing hypersurfaces). For a clear exposition of this theory we refer to the book of Giusti~\cite{Giu83}.
             
\begin{state}{\bf 7.3 Theorem.} %
  If \hskip1pt$G=\graph u$, where $u\in C^{2}(\R^{n}\setminus B_{1})$ satisfies the MSE on $\R^{n}\setminus
B_{1}$, then $G$ is asymptotically conic in the sense of~1.12(b); that is for each sequence $\rho_{k}\to \infty$
there is a subsequence $\rho_{k^{\prime}}$ and a cone $\C\in{\cal{}P}$ such that $\eta_{0,\rho_{k^{\prime}}}G\to
{\cal{}P}$ in $\R^{n+1}\setminus\{0\}$ in the sense of~1.4. In fact in this case we have always that $\C$ is
cylindrical: $\C=\C_{0}\times \R$ for some $(n-1)$-dimensional area minimizing cone $\C_{0}\subset \R^{n}$.
\end{state} %

{\bf{}Proof:} For the proof of this in the case when $u$ is an entire solution (i.e.\ when $u$ is defined over all of
$\R^{n}$), see for example~\cite{Mir77} and~\cite{Giu83}; the proof for exterior solutions requires only a little
more argument, and is described in~\cite{Sim87a}.

\medskip

Next we recall that if $\nu=(\nu^{1}\ldots, \nu^{n+1})={(-Du,1)\over \sqrt{1+|Du|^{2}}}$ is the upward pointing
unit normal of $G$ (thought of as a function on $G$ rather a function in $\R^{n}\setminus B_{1}$), then
$\nu^{n+1}\equiv e_{n+1}\cdot \nu$ satisfies the Jacobi-field equation
$$%   
\Delta_{G} u + q_{G}u =0,  
\leqno{7.4} 
$$% 
where $q_{G}=|A_{G}|^{2}$ is the square length of the second fundamental form of $G$.  Of course by the
$C^{\infty}$ convergence of $G_{k^{\prime}}$ to $\C=\C_{0}\times \R$, we trivially have that
$\rho_{k^{\prime}}^{2}q_{G_{k^{\prime}}}\to q_{\C}$ uniformly on compact subsets of $\C$, where $q_{\C}=
|A_{\C}|^{2}\equiv |A_{\C_{0}}|^{2}$ is the square length of the second fundamental form of $\C$, or,
equivalently, $\C_{0}$.  Thus all the conditions for the application of the Theorem~2 of~\S1 do hold (with
$a=0,b=0$ and $q=|A_{G}|^{2}$ in this case), and hence we conclude that
$$% 
R^{-n}\|u\|_{L^{1}(G\cap B_{\R}\setminus B_{R/2})} \le R^{-\gamma}
$$% 
for all sufficiently large $R$ where $\gamma$ is any number less than $\gamma_{0}$, where $\gamma_{0}$
is $\inf \gamma(\Sigma)$, where $\gamma(\Sigma)={n-2\over 2}-\sqrt{({n-2\over 2})^{2}+\lambda_{1}(\Sigma})$,
and the inf is over all $\Sigma=\C\cap S^{n}$ corresponding to all possible tangent cylinders $\C=\C_{0}\times
\R$ of $G$ at $\infty$; as in~1.8, $\lambda_{1}(\Sigma)$ is the first eigenvalue of the operator
$-(\Delta_{\Sigma}+ q_{\Sigma})$\renewcommand{\thefootnote}{\fnsymbol{footnote}}\footnote[3]{Actually
in~1.8 we defined $\lambda_{1}(\Sigma)$ to be maximum of $\lambda_{1}(\Sigma_{\ast})$ over all the (finitely
many) connected components $\Sigma_{\ast}$ of $\Sigma$, but in this case by a result of Bombieri and
Giusti~\cite{BomG72} we automatically have that $\Sigma$ is connected.}  in case $q_{\Sigma}$ is the square
length of the second fundamental form of $\Sigma$.

We actually claim that the exponent $\gamma_{0}$ above can be computed in terms of the first eigenvalues of
the cross-sectional cones as follows:

\begin{state}{\bf 7.5\, Lemma.} % 
$\gamma(\Sigma)= \gamma(\Sigma_{0})$, where $\Sigma=(\C_{0}\times
\R)\cap S^{n}$, $\Sigma_{0}=\C_{0}\cap S^{n-1}$, $\gamma(\Sigma)= {n-2\over 2}-\sqrt{({n-2\over
2})^{2}+\lambda_{1}(\Sigma})$ and $\gamma(\Sigma_{0})= {n-3\over 2}-\sqrt{({n-3\over
2})^{2}+\lambda_{1}(\Sigma_{0}})$.
\end{state}

{\bf{}Proof:} For $\delta>0$, select a $C^{\infty}$ relatively open $\Omega_{\delta}\subset\Sigma_{0}$ with
$\{\omega\in\Sigma_{0}:\dist(\omega,\sing\Sigma_{0})>\delta\}\subset \Omega_{\delta} \subset
\overline\Omega_{\delta}\subsetneq\Sigma_{0}$ and let $\varphi_{\delta}$ be the (smooth and positive) first
eigenfunction for the operator $-(\Delta_{\Sigma_{0}}+q_{\Sigma_{0}})$ with zero Dirichlet data, i.e.,
$\varphi_{\delta}=0$ on $\partial\Omega_{\delta}$. Then $\lambda_{1}(\Omega_{\delta})\to
\lambda_{1}(\Sigma_{0})$ as $\delta\downarrow 0$ and
$\lambda_{1}(\Omega_{\delta})>\lambda_{1}(\Sigma_{0})$ for every $\delta>0$ and so for $\delta$ sufficiently
small we have $-({n-3\over{}2})^{2}<\lambda_{\delta} <0$, and in particular this means that
$\gamma^{2}-(n-3)\gamma-\lambda_{\delta}$ has real roots, with smaller root $\gamma_{\delta}$,
$\gamma_{\delta}={n-3\over{}2}-\sqrt{({n-3\over{}2})^{2}+\lambda_{\delta}}$, satisfying
$$%
\gamma_{\delta} < \tfrac{n-3}{2} \text{ and }  %
   \gamma_{\delta}^{2}-(n-3)\gamma_{\delta}=\lambda_{1}(\Omega_{\delta}).  %
\leqno{(1)}
$$%
For $(x,y)\in \R^{n}\times\R=\R^{n+1}$, let $r_{0}=|x|$ and $r=\sqrt{|x|^{2}+y^{2}}$, and let
$$%
\Lambda_{\delta}= (\C_{0,\delta}\times \R)\cap  S^{n},
$$%
where
$$%
\C_{0,\delta}=\{r\omega:r>0, \omega\in\Omega_{\delta}\}\subset\C_{0}.
$$% 
Since
$\Delta_{\C_{0}}=r_{0}^{2-n}{\partial\over\partial r_{0}}\bigl(r_{0}^{n-2}{\partial \over
  \partial r_{0}}\bigr) + {1\over r_{0}^{2}}\Delta_{\Sigma_{0}}$, 
by direct computation we see that
$v_{\delta}\equiv r_{0}^{-\gamma_{\delta}}\varphi_{\delta}$ satisfies the equation
$$%
\Delta_{\C_{0}}v_{\delta}+ r_{0}^{-2}q_{\Sigma_{0}}v_{\delta}=0 
$$%
and hence of course it also is a solution (independent of the $y$-variable) of the equation $\Delta_{\C}v_{\delta}+
r_{0}^{-2}q_{\Sigma_{0}}v_{\delta}=0$. Indeed since $v_{\delta}$ can be written $r^{-\gamma_{\delta}}
(r/r_{0})^{\gamma_{\delta}}\varphi_{\delta}$ and since the Laplacian on $\C$ has the form
$r^{1-n}{\partial\over\partial r}\bigl(r^{n-1}{\partial \over
  \partial r}\bigr) + {1\over r^{2}}\Delta_{\Sigma}$, we see that $\Phi_{\delta}\equiv
(r/r_{0})^{\gamma_{\delta}}\varphi_{\delta}$ is a solution of the equation
$$%
-(\Delta_{\Sigma}+q_{\Sigma})\Phi_{\delta} = (\gamma_{\delta}^{2} -   %
                                          (n-2)\gamma_{\delta}) \Phi_{\delta}  \text{ on $\Lambda_{\delta}$,} %
\leqno{(2)}
$$%
where $q_{\Sigma}=r_{0}^{-2}q_{\Sigma_{0}}$.  Also since $\varphi_{\delta}$ is $C^{1}$ and
$\gamma_{\delta}<{n-3\over{}2}$ it is readily checked that $\Phi_{\delta}\in W_{0}^{1,2}(\Lambda_{\delta})$. 
Then we can use a standard argument to assert that
$$%
(\gamma_{\delta}^{2}-(n-2)\gamma_{\delta})=\lambda_{1}(\Lambda_{\delta}) %
\leqno{(3)}
$$%
where
$$%
\lambda_{1}(\Lambda_{\delta})= \inf_{v\in C^{1}_{\text{c}}(\Lambda_{\delta})\setminus\{0\}}  %
     \|v\|_{L^{2}(\Lambda_{\delta})}^{-2}   \int_{\Lambda_{\delta}}(|\nabla_{\Sigma}v|-q_{\Sigma}v^{2}),  %
$$%
as follows. First note that because $\Phi_{\delta}\in W_{0}^{1,2}(\Lambda_{\delta})$ we can multiply by
$\Phi_{\delta}$ in~(1) and integrate by parts, thus showing that
$$%
\gamma_{\delta}^{2} -   (n-2)\gamma_{\delta}  \ge \lambda_{1}(\Lambda_{\delta}).  %
\leqno{(4)}
$$%
Take any smooth subdomain $\Omega\subset\subset\Lambda_{\delta}$ and let $\lambda_{1}(\Omega)$ be the
minimum eigenvalue for $-(\Delta_{\Sigma}+q_{\Sigma})$ on $\Omega$. We claim that then
$\lambda_{1}(\Omega)>\lambda=\gamma_{\delta}^{2} - (n-2)\gamma_{\delta}$, because if
$\lambda_{1}(\Omega)\le \lambda$ then, with $\varphi\in C^{1}_{0}(\overline \Omega)$ the positive smooth
eigenfunction corresponding to the eigenvalue $\lambda_{1}(\Omega)$, we would have
$$%
\Delta_{\Sigma}(\Phi_{\delta}-\mu \varphi)+(q_{\Sigma}+\lambda)(\Phi_{\delta}-\mu \varphi)\le 0  %
                                                             \text{ on }\Omega  %
$$%
for $\mu >0$, and we could take $\mu>0$ such that $\Phi_{\delta}-\mu \varphi$ has a zero minimum in
$\Omega$, contradicting the maximum principle.  Thus we must have $\lambda_{1}(\Omega)>\lambda$ for all
such $\Omega$. Since $\lambda_{1}(\Lambda_{\delta})=\inf\lambda_{1}(\Omega)$ over all such $\Omega$ we
thus have $\lambda_{1}(\Lambda_{\delta})\ge \lambda$ and hence (by~(4)) we must have~(3) as claimed.

Then, by~(1) and~(3), $\gamma_{\delta}$ is the (smaller) root of both the equation
$\gamma^{2}-(n-2)\gamma-\lambda_{1}(\Lambda_{\delta})=0$ and also the equation
$\gamma^{2}-(n-3)\gamma-\lambda_{1}(\Omega_{\delta})=0$, so~7.5 follows by letting $\delta\downarrow 0$.

\bigskip

Finally we want to establish the bound $\lambda_{1}(\Sigma_{0})\le -(n-2)$ mentioned above:

\begin{state}{\bf 7.6 Lemma.} %
  With $\C=\C_{0}\times\R$ any tangent cone at $\infty$ for $G$ and with $\Sigma_{0}=\C_{0}\cap S^{n-1}$, we
have $\lambda_{1}(\Sigma_{0})\le -(n-2)$.
\end{state}

Recall that by the identity of James Simons \cite{Sis68} (see also \cite{SchSY75})
$\Delta_{\Sigma_{0}}|A_{\Sigma_{0}}|+ q_{\Sigma_{0}}|A_{\Sigma_{0}}|\ge (n-2)|A_{\Sigma_{0}}|$, and this
suggests that we should try to use~4.9.  (Notice that the coefficient of the term on right side is indeed $n-2$,
because $n-2$ is the dimension of $\Sigma_{0}$.)

To make it possible to apply~4.9 we need also to recall that by~\cite{SchSY75} we have the estimates
$$%   
\int_{\C_{0}}|A_{\C_{0}}|^{p}f^{p} \le C\int_{\C_{0}}|\nabla f|^{p} %
\leqno{(1)} %
$$% 
for $p\in [2,4+\sqrt{8/n})$, provided $f\in C^{\infty}_{c}(\C_{0})$.  (Notice that the work of~\cite{SchSY75}
formally assumed $\sing \C_{0}=\{0\}$, but the proof of course works without change in the case of an arbitrary
singular set, so long as we assume, as we do here, that $\support f$ is a compact subset of the smooth manifold
$\C_{0}$.)  Now recall that in the present area minimizing case (see e.g.~\cite{Giu83}, keeping in mind that the
dimension of $\C_{0}$ is $n-1$) we have that $\dim\sing \C_{0}\le n-8$, so that in particular
$$% 
{\cal{}H}^{n-5}(\sing\C_{0})=0, \; \text{ hence } {\cal{}H}^{n-6}(\sing\Sigma_{0})=0.  \leqno{(2)}
$$% 
In view of~(1) (with $p=4$) and~(2) we can use precisely the same argument as in~\S4 (preceding~4.9) in order
to deduce that~(1) is also valid for any $f\in C^{1}_{c}(\breve B_{R}\setminus B_{\tau})$ for any
$0<\tau<R<\infty$.  Hence we conclude from~(1) that $A_{\C_{0}}\in L^{4}(\C\cap (B_{2}\setminus B_{1/2}))$,
and hence $A_{\Sigma_{0}}\in{}L^{4}(\Sigma_{0})$. Therefore we can apply~4.9 to conclude that
$$%   
\lambda_{1}(\Sigma_{0}) \le -(n-2)
$$% 
as claimed.

\medskip

Thus with ${\cal{}E}_{0}=\{\C_{0}\cap S^{n-1}:\C_{0}\times\R\text{ is a tangent cylinder of $G$ at $\infty$}\}$,
we have
$$%
\lambda_{1}({\cal{}E}_{0})\le -(n-2)  
$$%
and hence by Lemma~7.5 we can apply the main decay estimate of Theorem~2 in~\S1 with
$\gamma_{0}={n-3\over{}2}-\sqrt{({n-3\over{}2})^{2}-(n-2)}$ to the solution $u=\nu_{n+1}$ (as in~7.4), whence
we conclude that for each $\gamma<\gamma_{0}$ there is a constant $\rho_{0}>1$ with
 $$
 \int_{S_{R}\setminus S_{R/2}} \nu_{n+1}\le R^{-\gamma},\quad R\ge \rho_{0}.
 $$
 On the other hand by the H\"older inequality and the standard lower bound on the volume of the intersection
of the graph $G$ with an $(n+1)$-dimensional ball we have
\begin{align*}   %
  CR^{n} &\le {\cal{}H}^{n}(S_{R}\setminus S_{R/2})\\ &\le \int_{S_{R}\setminus
S_{R/2}}\nu_{n+1}\,d{\cal{}H}^{n} \int_{S_{R}\setminus S_{R/2}} \sqrt{1+|Du|^{2}}\,d{\cal{}H}^{n}\\ &\le
CR^{-\gamma}\int_{S_{R}\setminus S_{R/2}}\sqrt{1+|Du|^{2}}\,d{\cal{}H}^{n},
\end{align*} %
whence
 $$
 R^{-n}\int_{S_{R}} |Du| \,d{\cal{}H}^{n} \ge C R^{\gamma},\quad R\ge \rho_{0},
 $$
 as claimed in the introduction.

\bigskip

\bibliographystyle{amsalpha}

\newcommand{\noopsort}[1]{}
\providecommand{\bysame}{\leavevmode\hbox to3em{\hrulefill}\thinspace}
\providecommand{\MR}{\relax\ifhmode\unskip\space\fi MR }
\providecommand{\MRhref}[2]{%
  \href{http://www.ams.org/mathscinet-getitem?mr=#1}{#2}
}
\providecommand{\href}[2]{#2}

\end{document}